\documentclass{amsart}
\newtheorem{theorem}{Theorem}[section]

\usepackage{graphicx,palatino,pifont,times}
\usepackage{epsfig}
\usepackage{a4wide}
\usepackage{subfig}
\theoremstyle{definition}
\newtheorem{definition}[theorem]{Definition}

\newtheorem{proposition}[theorem]{Proposition}

\theoremstyle{remark}
\newtheorem{remark}[theorem]{Remark}

\newtheorem{corollary}[theorem]{Corollary}

\numberwithin{equation}{section}

\begin{document}
	
	\title{A study on sensitivity and stability analysis of non-stationary $\alpha$-fractal functions}
	
%

	
\author{Anarul Islam Mondal}
\address{Department of Mathematics, National Institute of Technology Rourkela,	Rourkela 769008}
\email{anarulmath96@gmail.com}

\author{Sangita Jha}
\address{Department of Mathematics, National Institute of Technology Rourkela,	Rourkela 769008}
\email{jhasa@nitrkl.ac.in}

\begin{abstract}
\noindent  This article aims to study fractal interpolation functions corresponding to a sequence of iterated function systems (IFSs). For a suitable choice of a sequence of IFS parameters, the corresponding non-stationary fractal function is a better approximant for the non-smooth approximant. In this regard, we first construct the non-stationary interpolant in the Lipschitz space and study some topological properties of the associated non-linear fractal operator. Next, we discuss the stability of the interpolant having small perturbations. Also, we investigate the sensitivity with respect to the perturbations of the IFS parameters by providing an upper bound of errors acquired in the approximation process.
In the end, we study the continuous dependence of the proposed interpolant on different IFS parameters.
\end{abstract}

\maketitle

\section{Introduction}
The representation of arbitrary functions or data sets in terms of simple classical functions like polynomials, trigonometric functions, or exponentials is one of the main ideas of numerical analysis and approximation theory. Traditional approaches may not produce an approximant with the required precision when dealing with irregular forms, such as real-world signals like time series, financial series, climatic data, and bioelectric recordings.  To address the irregularity in different practical situations, Barnsley \cite{Barnsley} introduced the notion of fractal interpolation functions (FIFs). The reader can consult the books \cite{Barnsley1, Massopustbook} and the sources provided there for a clear explanation of several crucial subjects in this direction.

FIFs have some advantages over traditional interpolation functions. FIFs, in general, are self-similar/affine, and the Hausdorff Besicovitch dimensions of their graphs are non-integers. The reader can consult \cite{Nasim, SS, RuanSuYao, VM} for dimensional analysis of stationary fractal functions. The primary advantage of fractal interpolants over classical interpolants is generating both smooth and nonsmooth interpolants, depending on the choice of the IFS parameters. 
It is important to note that relatively few techniques-including subdivision schemes, another well-liked technique, produce nonsmooth interpolant \cite{DynLevin}. Subdivision methods and fractal interpolation have recently been linked in some attempts \cite{DynLevinMassopust, LevinDynViswanathan}.

Motivated by the work of Barnsley, a family of non-affine fractal functions $f^\alpha$-also known as $\alpha$-fractal functions corresponding to a continuous function $f$ on a closed and bounded interval $I$ of $\mathbb{R}$ was studied by Navascu\'es \cite{Navascues1}. The function $f^\alpha$ approximates and interpolates $f$.  This method of fractal perturbation produces an operator,  which links the theory of FIF to the area of Operator theory, Approximation theory, Functional analysis, and Harmonic analysis.

Evidently, a lot of work has been done in the development of fractal approximations (stationary) \cite{Vishal, BouboulisDalla, Vij, Zipper, Abbas, Ri, Secelean}, and it is still a subject of very active research, with an extensive list of connections and applications. But on the other hand, many questions remain to be settled in the non-stationary case, and more specifically in the case of a sequence of IFSs. The difficulty lies in the fact that the standard tools are not well developed in the setting of a sequence of IFSs. Consequently, advances in non-stationary IFS are required. 


Using the fixed points of contractive operators for a specific type of IFS is a helpful way to build fractals \cite{Hutchinson}. In the case of usual fractal functions, a common procedure of construction is to consider an operator on a complete space and define the fractal function as the unique fixed point \cite{ViswanathanNavascues}. However, if the operator is substituted by a sequence of maps, then we face problems with similar construction. 
Levin et al.\cite{LevinDynViswanathan} investigate the trajectory of contraction mappings which produces limit attractors at various scales with different features or shapes. Recently, Massopust \cite{Massopust} advanced fractal function to a new and more adaptive environment by utilising the concept of forward and backward trajectories to offer new forms of fractal functions with the variable local and global behaviours. 

In this article, we construct a non-stationary fractal interpolation function on the space of all Lipschitz functions on a closed and bounded interval $I$ of $\mathbb{R}$. Our goal in presenting the non-stationary variant is not only to generalize the stationary situation; instead, we hope to extend the applicability of non-stationary fractal functions. A loss of flexibility and visible inaccuracies in the fitting and approximation of some complex curves and non-stationary data that show less self-similarity may be caused by the obvious self-similarity characteristics of FIFs created by IFSs with constant parameters. The non-stationary FIF, in general, does not depend on local or global data points; we may use non-stationary FIFs to get more flexibility and accuracy. We look into the stability and sensitivity of the analytical characteristics of the FIFs generated by the class of IFSs with changeable parameters.
Further, we calculate an upper bound of the errors that occurred when we compare our proposed interpolant with the non-stationary interpolant that came from the perturbed IFS.

The IFS plays an important role in the theory of fractals. In particular, for application purposes, like in computer graphics, IFS with parameters are employed. Therefore, it is interesting to verify if the small changes in the parameters lead to small changes in the interpolant. 
In \cite{Barnsley1}, Barnsley proved that the attractor depends continuously on the parameters with respect to the Hausdorff metric if all the IFS maps depend continuously on the parameters. The fact is very important as it tells us to interpolate smoothly between atttractors which is helpful in image animation and computer graphics. In this regard, we study the continuous dependency of the IFS parameters of our proposed non-stationary interpolant, which also generalize the same study in the stationary case \cite{VermaViswanathan}. 

The article is organized as follows. We review the concepts of trajectories and stationary as well as non-stationary IFSs in Section 2. In Section 3, we describe the construction of the non-stationary $\alpha$-fractal functions on the Lipschitz space $Lip_d(I)$. A nonlinear fractal operator on $Lip_d(I)$ is explored, and some essential topological properties of the operator are discussed in Section 4. In Section 5, we study stability analysis by perturbing the ordinates. Also, we investigate the sensitivity with respect to the perturbations of the IFS parameters by providing an upper bound of errors acquired in the approximation process. In the end, we investigate the continuous dependence of non-stationary FIFs on IFS parameters $\alpha,$ $b$ ~\text{and}~ $\triangle.$

\section{Preparatory Results}
In this section, we collect the necessary tools leading to the construction of non-stationary $\alpha$-fractal functions in the Lipschitz space. For more details on this section, we invite the reader to study the paper of Massopust \cite{Massopust} and the work by Barnsley \cite{Barnsley} and Navascu\'es \cite{Navascues1}.
\subsection{Iterated Function System and Trajectories}
Let $(\mathbf{\mathbf{X}},d)$ be a complete metric space. Let $\mathcal{H}(\mathbf{\mathbf{X}})$ denote the collection of all non-empty compact subsets of $\mathbf{\mathbf{X}}$ and define the Hausdorff distance between sets $A$ and $B$ of $\mathcal{H}(\mathbf{\mathbf{X}})$ as 
$$d_{H}(A,B) = \inf \{ \epsilon \ge 0 ; A \subseteq B_\epsilon~\text{and}~B \subseteq A_\epsilon \},$$
where  $A_\epsilon =\displaystyle \bigcup_{x \in A} \{ z \in \mathbf{X} ;~ d(z,x) \le \epsilon \}$.
The space $(\mathcal{H}(\mathbf{\mathbf{X}}), d_H)$ is a complete metric space known as the space of fractals.

\begin{definition}
Let $(\mathbf{\mathbf{X}},d)$ be a complete metric space and $ w_{i}:\mathbf{X} \to \mathbf{X}$ be $N$ continuous maps. Then the system $\mathcal{I} = \{ \mathbf{X} ; w_{i} : i= 1,2,\dots,N\}$ is called an iterated function system (IFS).
If each map $w_{i}$ in $\mathcal{I}$ is a contraction, then the IFS $\mathcal{I}$ is hyperbolic.
\end{definition}
The attractor of the IFS is the fixed point of the set-valued Hutchinson map $\mathcal{W}: \mathcal{H}(\mathbf{X}) \to \mathcal{H}(\mathbf{X})$,
$$\mathcal{W}(E) = \bigcup_{i=1}^{N} w_{i}(E),~~~~~~ E \in \mathcal{H}(\mathbf{X}).$$ 
$W$ is a contraction map on $\mathcal{H}(\mathbf{X})$ with the Lipschitz constant
$Lip(\mathcal{W}):= \max \{ Lip(w_{i}):i=1,2,\dots,N\}.$ The uniqueness of the fixed point is guaranteed by the Banach fixed point theorem. The limit of the iterative process $A_{k} = \mathcal{W}(A_{k-1}); k \in \mathbb{N}$, where $A_{0} \in \mathcal{H}(\mathbf{X})$ is any arbitrary set, also gives the attractor of the IFS.
\begin{definition}
    Let $T: \mathbf{X} \to \mathbf{X}$ be a contraction map on a complete metric space $(\mathbf{X},d)$. The forward iterates of $T$ are transformations $T^{\circ\,n} : \mathbf{X} \to \mathbf{X}$ defined by $T^{\circ\,0} (x) = x$ and $T^{\circ\,(n+1)}(x) = T(T^{\circ\,n}(x)) = \underbrace{T\circ T\circ \dots \circ T}_{( n+1~ \text{times})}(x)$~~ for $n \in \mathbb{N}\cup\{0\}$.
\end{definition}
\begin{definition}(Forward and Backward Trajectories)
Let $\mathbf{X}$ be a metric space and
$\{T_{r}\}_{r \in \mathbb{N}}$ be a sequence of Lipschitz maps on $\mathbf{X}$. We define forward and backward trajectories respectively
$$ \phi_r := T_{r}\ \circ \ T_{r-1}\ \circ\ \dots\ \circ\ T_{1} \ \ \text{and} \ \ \psi_r := T_{1}\ \circ\ T_{2}\ \circ\ \dots\ \circ\ T_{r}.$$
\end{definition}
The subject now concerns the convergence of general trajectories, i.e., under what conditions the forward and backward trajectories will converge. In addition, we look for the trajectories producing new types of fractal sets.
Recently, Levin et al. \cite{LevinDynViswanathan} observed that the backward trajectories converge under relatively mild conditions and may produce a new class of fractal sets. In \cite{LevinDynViswanathan}, the authors use the assumption of compact invariant domain to guarantee the convergence of backward trajectories. In \cite{NavascuesVerma}, Navascu\'es and Verma replace the compact invariant domain condition by a weaker condition. We now recall the result in the following.


\begin{proposition}\cite[Proposition 2.6]{NavascuesVerma} \label{Thm-nonstationary}
Let $\{T_r\}_{r \in \mathbb{N}}$ be a sequence of Lipschitz maps on a complete metric space $(\mathbf{X},d)$ with Lipschitz constant $c_r$. If~ $\exists~ x^{*}$~ in the space such that the sequence $\{ d( x^* , T_{r}(x^*))\}$ is bounded, and $\displaystyle \sum_{r=1}^{\infty} \prod_{i=1}^{r} c_i < \infty ,$ then the sequence $\{ \psi_{r}(x) \}$ converges for all $x \in \mathbf{X}$ to a unique limit $\bar{x}$.
\end{proposition}

\subsection{Stationary Fractal Interpolation Functions}
For a fixed $ k \in \mathbb{N}$, we denote by $\mathbb{N}_{k}$ the first $k$ natural numbers and $\mathbb{N}_{k}^{0} = \mathbb{N}_{k} \cup \{0\} $.
Let $I=[a,b]$ and define a partition $\Delta$ on $I$ by 
$$ \Delta = \{ ( x_0, x_1, \dots, x_N ) : a=x_0 < x_1 <\dots< x_N=b \}.$$
For $i \in \mathbb{N}_N$, let $I_i = [x_{i-1},x_{i}]$.
Suppose the affine maps $l_{i} : I \to I_{i}$ are such that 

\begin{equation}\label{stationary affine condition}
    l_{i}(x_0)=x_{i-1}~,~l_{i}(x_N)=x_{i},~ \text{and}~ |l_{i}(z_1) - l_{i}(z_2)| \le l |z_1 - z_2|~\forall z_1,\ z_2 \in I,
\end{equation}
where $0 \le l < 1$. Set $\mathcal{M} = I \times \mathbb{R}$. Let the $N$ continuous maps $F_{i} : \mathcal{M} \to \mathbb{R},~ i \in \mathbb{N}_N$ be such that 
\begin{equation}\label{stationary F-cont condition}
    F_{i}(x_0,y_0) = y_{i-1},\ \ \ \ \  F_{i}(x_N,y_N) = y_{i},
\end{equation}
and $F_{i}$ is a contraction in the second variable for $i \in \mathbb{N}_N.$ The following are the most chosen maps for the formation of IFS
\begin{equation}\label{stationary IFS-eq1}
    \begin{cases}
l_{i}(x) = a_{i}x+ e_{i},\\
F_{i}(x,y) = \alpha_{i}(x) y + q_{i}(x),\ \ \  i \in \mathbb{N}_N,
\end{cases}
\end{equation}
where $a_{i}, e_{i}$ can be determined by conditions (\ref{stationary affine condition}). $\alpha_{i}(x)$ are scaling functions satisfying $\|\alpha_i\|_\infty < 1$ and $q_{i}(x)$ are suitable continuous functions such that the condition (\ref{stationary F-cont condition}) is satisfied.
For $i \in \mathbb{N}_{N}$, we define  $W_{i} : \mathcal{M} \to I_{i} \times \mathbb{R}$ by
$$W_{i}(x,y) = (l_{i}(x) , F_{i}(x,y))\ \ \ \forall (x,y) \in \mathcal{M}.$$
The complete metric space $(\mathcal{M},d_H)$ with the above $N$-maps $\{W_i : i \in \mathbb{N}_N\}$ forms an IFS. The uniqueness of the attractor of the IFS is given by Barnsley; mentioned below.
\begin{theorem} \cite{Barnsley} \label{Barnsley thm}
    The IFS $\mathcal{I} = \{ \mathcal{M}; W_{i}:i \in \mathbb{N}_{N}\}$ admits a unique attractor $G$ and $G$ is the graph of a continuous function $f: I \to \mathbb{R}$ which interpolates the given data.
\end{theorem}
\subsection{Stationary $\alpha$-fractal Functions}
Let $\mathcal{C}(I)$ be the set of all continuous functions on a compact interval $I$. Let $f \in \mathcal{C}(I)$ be a prescribed function. Let us choose the function $q_i(x)$ in (\ref{stationary IFS-eq1}) as 
\begin{equation*}
    q_i(x) = f(l_{i}(x))-\alpha_{i}(x) b(x), ~ i \in \mathbb{N}_N,
\end{equation*}
where $\alpha_{i}: I \to \mathbb{R}$ are continuous scaling functions satisfying $\|\alpha_i\|_\infty < 1$ and $b: I \to \mathbb{R}$ be continuous function such that $b \neq f$, $b(x_0) = f(x_0)$ and $b(x_N) = f(x_N)$.
By Theorem \ref{Barnsley thm}, the corresponding IFS has a unique attractor $G$, which is the graph of a continuous map, namely $f_{\Delta,b}^{\alpha}: I \to \mathbb{R}$ that interpolates the given data set.
The map is known as the $\alpha$-fractal function.

Though the stationary fractal interpolation interpolates irregular functions very well, it depends on local data points. To get fractal interpolation that is independent of local data points and gets more reflexibility, Massopust \cite{Massopust} introduced non-stationary fractal interpolation by taking a sequence of IFSs. In \cite{NavascuesVerma}, the authors studied the parameterized non-stationary FIFs in $\mathcal{C}(I)$.

\subsection{Non-stationary $\alpha$-fractal Functions}\label{sec-2.4}
Let $f \in \mathcal{C}(I)$ and $r \in \mathbb{N}$. We use the following notation:
\begin{center}
$\alpha_{r}:= ( \alpha_{1,r} , \alpha_{2,r},\dots,\alpha_{N,r}) ,\ \  \alpha := \{ \alpha_r \}_{r \in \mathbb{N}}$\ \ and\ \ $b:= \{ b_r \}_{r \in \mathbb{N}}.$
\end{center}
Let $\alpha_{i,r} :I \to \mathbb{R}$ be continuous functions such that $$\|\alpha\|_\infty = \sup\{\|\alpha_{r}\|_\infty : r \in \mathbb{N}\} < 1, \ \ \text{where}\ \ \|\alpha_{r}\|_\infty = \sup \{\|\alpha_{i,r}\|_\infty : i \in \mathbb{N}_N\},$$
 and $b_r \in \mathcal{C}(I)$ such that 
 \begin{equation} \label{base cond}
     b_r \neq f,~~~b_r(x_0)=f(x_0)~~~ \text{and}~~~b_r(x_N)= f(x_N).
 \end{equation}
To define a sequence of IFSs, we use the following sequence of continuous maps 
\begin{equation}\label{non-stationary IFS-eq1}
    \begin{cases}
l_{i}(x) = a_{i}x+ e_{i} = \frac{x_i-x_{i-1}}{x_N-x_0}x+\frac{x_Nx_{i-1}-x_0x_i}{x_N-x_0},\\
F_{i,r}(x,y) = \alpha_{i,r}(x) y + f(l_{i}(x))-\alpha_{i,r}(x) b_r(x),\ \ \  i \in \mathbb{N}_N.
\end{cases}
\end{equation}
For each $i \in \mathbb{N}_{N}$, we define $$W_{i,r} : \mathcal{M} \to I_{i} \times \mathbb{R} \;  \text{by} \;  W_{i,r}(x,y) = (l_{i}(x) , F_{i,r}(x,y)).$$
Now we have a sequence of IFSs
$\mathcal{I}_r = \{ \mathcal{M};W_{i,r} : i \in \mathbb{N}_{N}\} .$
\noindent 
Let $$ \mathcal{C}_{f}(I) = \{ g \in \mathcal{C}(I) : g(x_i)=f(x_i),\ i= 0,N \}.$$
Then $\mathcal{C}_{f}(I)$ is a complete metric space.
For $r \in \mathbb{N}$, we define a sequence of Read-Bajraktarevi\'c (RB) operators $ T^{\alpha_r} : \mathcal{C}_{f}(I) \to \mathcal{C}_{f}(I)$ by 
\begin{equation} \label{RB operator eq}
    \begin{aligned}
       (T^{\alpha_r} g)(x) =&~ F_{i,r}( {l_{i}}^{-1}(x) ,~g({l_{i}}^{-1}(x)))\\ =&~ f(x) + \alpha_{i,r}(Q_{i}(x))\cdot g(Q_{i}(x))- \alpha_{i,r}(Q_{i}(x))\cdot b_{r}(Q_{i}(x)),
    \end{aligned}
\end{equation}
for $x \in I_{i},~~i \in \mathbb{N}_N$, where $Q_{i}(x) = l_{i}^{-1}(x)$.\\
The above operator is well defined and for any function $h \in \mathcal{C}_f(I)$, the sequence of backward trajectories $\{ T^{\alpha_1}\ o\ T^{\alpha_2}\ o \dots o\ T^{\alpha_r} h \}$ converges to a map $f^{\alpha}$ of $C_f (I)$ \cite{NavascuesVerma}. The map $f^{\alpha}$ is the unique map that satisfies the following equation
\begin{equation} \label{self referential eq}
    f^{\alpha} (x) = f(x) + \lim_{r \rightarrow \infty} \sum_{j=1}^{r} \alpha_{i,1}(Q_{i}(x))\dots\alpha_{i,j}(Q_{i}^{j}(x)) (f-b_j)(Q_{i}^{j}(x)).
\end{equation}
$f^{\alpha}$ is called the non-stationary $\alpha$-fractal interpolation function.


\section{Non-stationary $\alpha$-fractal function on Lipschitz Space}
\noindent Let $g : I \to \mathbb{R}$ be a function. For $0 < d \le 1$, define \[Lip_d(g)\ =\ \sup \left\{\dfrac{|g(x) - g(y)|}{|x-y|^d} : x,y \in I\ \text{and}\ x \neq y \right\}.\]
The Lipschitz space is defined as $Lip_d(I)\ =\ \left\{ g : I \to \mathbb{R}\ :\ \ Lip_d(g) < \infty \right\}.$
Define $\|g\|_{d} = \max \{ \|g\|_\infty , Lip_d(g) \}$. It is routine to show that $(Lip_d(I) , \|.\|_{d})$ is a Banach space. For more details of the Lipschitz functions in an arbitrary Banach space, please refer to \cite{Johnson}. 
Let \[Lip_{d,f}(I) =\{ g \in Lip_d(I) : g(x_0) = f(x_0), g(x_N) = f(x_N) \} .\]
Then being a closed subspace of the Banach space $Lip_{d}(I)$, $Lip_{d,f}(I)$ is a Banach space.
\begin{theorem}
Let $f \in Lip_d(I)$. Let $r \in \mathbb{N}, b_r \in Lip_{d,f}(I)$ be such that $\|b\|_{d}: =\sup\limits_{r \in \mathbb{N}} \|b_r\|_{d} < \infty$ and the scaling functions $\alpha_{i,r}\in Lip_{d}(I)$ are chosen such that $\max\limits_{i \in \mathbb{N}_{N}} \left(\frac{\|\alpha_{i,r}\|_{d}}{a_{i}^{d}}\right) < \dfrac{1}{2}$. 
We define a sequence of RB operators $T^{\alpha_r} : Lip_{d,f}(I) \to Lip_{d,f}(I)$ by 
\begin{equation} \label{Lip RB operator eq}
       (T^{\alpha_r} g)(x) =f(x) + \alpha_{i,r}(Q_{i}(x))\cdot g(Q_{i}(x))- \alpha_{i,r}(Q_{i}(x))\cdot b_{r}(Q_{i}(x)),
\end{equation}
for $x \in I_{i},~~i \in \mathbb{N}_N$.
Then the following hold.
\begin{enumerate}
\item The RB operator $T^{\alpha_r}$ defined in equation (\ref{Lip RB operator eq}) is well defined on $Lip_{d,f}(I)$.
\item In fact, $T^{\alpha_r} : Lip_{d,f}(I) \to Lip_{d,f}(I) \subset Lip_{d}(I)$ is a contraction map.
\item There exists a unique function $f_{b,Lip_d}^{\alpha} \in Lip_{d,f}(I)$ such that the sequence $\{ T^{\alpha_1}\ o\ T^{\alpha_2}\ o \dots o\ T^{\alpha_r} g \}$ converges to the map $f_{b,Lip_d}^{\alpha}$ for every $g \in Lip_{d,f}(I)$.
\end{enumerate}
\end{theorem}
\begin{proof}
\begin{enumerate}
    \item The norm defined on $Lip_{d,f}(I)$ is $\|f\|_{d} = \max \{ \|f\|_\infty , Lip_{d}(f) \}$ so that, $ \|T^{\alpha_r} f\|_{d} = \max \{ \|T^{\alpha_r} f\|_\infty , Lip_{d}(T^{\alpha_r} f) \}.$ From the definition of RB operators, we have
 \[(T^{\alpha_r} g)(x) = f(x) + \alpha_{i,r}(Q_{i}(x))\cdot (g - b_{r})(Q_{i}(x)).\]  
Now, 
\allowdisplaybreaks
\begin{align*}
    Lip_{d}(T^{\alpha_r} g)&= \displaystyle \sup_{\substack{x,y \in I \\x \neq y }} \dfrac{|T^{\alpha_r} g(x) - T^{\alpha_r} g(y) |}{|x-y|^{d}}\nonumber\\
    &= \displaystyle \sup_{\substack{x,y \in I_{i} \\x \neq y }} \dfrac{ |f(x)-f(y) + \alpha_{i,r}(Q_{i}(x)). (g - b_{r})(Q_{i}(x)) - \alpha_{i,r}(Q_{i}(y)). (g - b_{r})(Q_{i}(y))|}{|x-y|^{d}}\nonumber\\
    &\le \displaystyle \sup_{\substack{x,y \in I \\x \neq y }} \dfrac{|f(x)-f(y)|}{|x-y|^{d}} +\max_{i \in \mathbb{N}_{N}} \displaystyle \sup_{\substack{x,y \in I_{i} \\x \neq y }} \dfrac{ |\alpha_{i,r}(Q_{i}(x)). [(g - b_{r})(Q_{i}(x)) - (g - b_{r})(Q_{i}(y))]|}{|x-y|^{d}}\nonumber\\
    & \hspace{3.33cm}  +\max_{i \in \mathbb{N}_{N}} \displaystyle \sup_{\substack{x,y \in I_{i} \\x \neq y }} \dfrac{|(g - b_{r})(Q_{i}(y))[\alpha_{i,r}(Q_{i}(x))-\alpha_{i,r}(Q_{i}(y))]|}{|x-y|^{d}}\nonumber\\
    &\le Lip_{d}(f) + \max_{i \in \mathbb{N}_{N}} \|\alpha_{i,r}\|_{\infty} \displaystyle \sup_{\substack{x,y \in I_{i} \\x \neq y }} \dfrac{|(g(Q_{i}(x))-g(Q_{i}(y)))| + |(b_{r}(Q_{i}(x)) - b_{r}(Q_{i}(y)))|}{|x-y|^{d}}\nonumber\\
    & \hspace{1.65cm}+ \|g-b_r\|_\infty \max_{i \in \mathbb{N}_{N}} \displaystyle \sup_{\substack{x,y \in I_{i} \\x \neq y }} \dfrac{|(\alpha_{i,r}(Q_{i}(x)) - \alpha_{i,r}(Q_{i}(y)))|}{|x-y|^{d}}\nonumber\\
    &= Lip_{d}(f) + \max_{i \in \mathbb{N}_{N}} \|\alpha_{i,r}\|_\infty \displaystyle \sup_{\substack{x,y \in I_{i} \\x \neq y }} \dfrac{|(g(Q_{i}(x))-g(Q_{i}(y)))| + |(b_{r}(Q_{i}(x)) - b_{r}(Q_{i}(y)))|}{a_{i}^{d}|Q_{i}(x)-Q_{i}(y)|^{d}}\nonumber\\
    & \hspace{1.65cm}+ \|g-b_r\|_\infty \max_{i \in \mathbb{N}_{N}} \displaystyle \sup_{\substack{x,y \in I_{i} \\x \neq y }} \dfrac{|(\alpha_{i,r}(Q_{i}(x)) - \alpha_{i,r}(Q_{i}(y)))|}{a_{i}^{d}|Q_{i}(x)-Q_{i}(y)|^{d}}\nonumber\\
    &=Lip_{d}(f) + \max_{i \in \mathbb{N}_{N}} \left(\frac{\|\alpha_{i,r}\|_\infty}{a_{i}^{d}}\right) \displaystyle \sup_{\substack{\Tilde{x},\Tilde{y} \in I \\ \Tilde{x} \neq \Tilde{y} }} \dfrac{|g(\Tilde{x})-g(\Tilde{y})| + |b_{r}(\Tilde{x}) - b_{r}(\Tilde{y}))|}{|\Tilde{x}-\Tilde{y}|^{d}}\nonumber\\
    & \hspace{1.65cm}+ \|g-b_r\|_\infty \max_{i \in \mathbb{N}_{N}} \displaystyle \sup_{\substack{\Tilde{x},\Tilde{y} \in I \\ \Tilde{x} \neq \Tilde{y} }} \dfrac{|(\alpha_{i,r}(\Tilde{x}) - \alpha_{i,r}(\Tilde{y}))|}{a_{i}^{d}|\Tilde{x}-\Tilde{y}|^{d}}\nonumber\\
    &\leq Lip_{d}(f) + \max_{i \in \mathbb{N}_{N}} \left(\frac{\|\alpha_{i,r}\|_\infty}{a_{i}^{d}}\right) ( Lip_{d}(g) + Lip_{d}(b_r) ) + \max_{i \in \mathbb{N}_{N}} \left(\frac{Lip_{d}(\alpha_{i,r})}{a_{i}^{d}}\right) (\|g\|_\infty+\|b_r\|_\infty)\nonumber\\
    &\le Lip_{d}(f) + \max_{i \in \mathbb{N}_{N}} \left(\frac{\|\alpha_{i,r}\|_{d}}{a_{i}^{d}}\right) ( Lip_{d}(g) + Lip_{d}(b_r) + \|g\|_\infty+\|b_r\|_\infty).
\end{align*}
As $f, g, b_r \in Lip_{d,f}(I)$, the above estimation ensures that $Lip_{d}(T^{\alpha_r} g) <\infty$ and so that $T^{\alpha_r} g \in Lip_{d}(I)$. Also $T^{\alpha_r} g(x_0) = f(x_0)$ and $T^{\alpha_r} g(x_N) = f(x_N)$. Hence $T^{\alpha_r} g \in Lip_{d,f}(I)$ and the RB operator $T^{\alpha_r}$ defined in equation (\ref{RB operator eq}) is well defined on $Lip_{d,f}(I)$.
    \item For $ x \in I_{i} $,
    \begin{equation*}
    \begin{split}
    |(T^{\alpha_r}g_1 - T^{\alpha_r}g_2)(x)| &= |\alpha_{i,r}(Q_i(x))\|(g_1 - g_2)(Q_i(x))|\\
    &\le \max_{i \in \mathbb{N}_{N}} (\|\alpha_{i,r}\|_\infty) \|g_1 - g_2\|_\infty,
    \end{split}
    \end{equation*}
    and hence
    \begin{equation} \label{norm difference}
        \|(T^{\alpha_r}g_1 - T^{\alpha_r}g_2)\|_\infty \le \max_{i \in \mathbb{N}_{N}} (\|\alpha_{i,r}\|_\infty) \|g_1 - g_2\|_\infty .
    \end{equation}
    Using similar steps in the estimation of $Lip_{d}(T^{\alpha_r} g)$, we obtain
    \begin{equation} \label{lip difference}
        Lip_{d}(T^{\alpha_r}g_1 - T^{\alpha_r}g_2) \le \max_{i \in \mathbb{N}_{N}} \left(\frac{\|\alpha_{i,r}\|_{d}}{a_{i}^{d}}\right) (Lip_{d}(g_1 - g_2) + \|g_1 - g_2\|_\infty).
    \end{equation}
    Combining (\ref{norm difference}) and (\ref{lip difference}), we get
    \begin{align*}
        &\|(T^{\alpha_r}g_1 - T^{\alpha_r}g_2)\|_{d}\\ &= \max \left\{\|(T^{\alpha_r}g_1 - T^{\alpha_r}g_2)\|_\infty , Lip_{d}(T^{\alpha_r}g_1 - T^{\alpha_r}g_2)\right\}\\
        &\le \max \left\{\max_{i \in \mathbb{N}_{N}} (\|\alpha_{i,r}\|_\infty) \|g_1 - g_2\|_\infty , \max_{i \in \mathbb{N}_{N}} \left(\frac{\|\alpha_{i,r}\|_{d}}{a_{i}^{d}}\right) (Lip_{d}(g_1 - g_2) + \|g_1 - g_2\|_\infty)\right\}\\
        &\le \max_{i \in \mathbb{N}_{N}} \left(\frac{\|\alpha_{i,r}\|_{d}}{a_{i}^{d}}\right) \max \left\{ \|g_1 - g_2\|_{d} ,  2\|g_1 - g_2\|_{d} \right\}\\
        &= 2 \max_{i \in \mathbb{N}_{N}} \left(\frac{\|\alpha_{i,r}\|_{d}}{a_{i}^{d}}\right) \|g_1 - g_2\|_{d}.
    \end{align*}
    By assumptions on the sequence of scaling functions, we can ensure that $T^{\alpha_r}$ is a contraction.

    \item Let $g \in Lip_{d}(I)$ be an arbitrary function. We have to check if the sequence $\{\| T^{\alpha_r}g - g\|_{d} \}$ is bounded.
    Now,
    \begin{align*}
    Lip_{d}(T^{\alpha_r} g - g) &= \displaystyle \sup_{\substack{x,y \in I \\x \neq y }} \dfrac{|(T^{\alpha_r} g - g )(x) - (T^{\alpha_r} g - g )(y) |}{|x-y|^{d}}\\
    &\le \displaystyle \sup_{\substack{x,y \in I \\x \neq y }} \dfrac{|T^{\alpha_r} g(x) - T^{\alpha_r} g(y) | + | g(x) - g(y)| }{|x-y|^{d}}\\
    &= \displaystyle \sup_{\substack{x,y \in I \\x \neq y }} \dfrac{|T^{\alpha_r} g(x) - T^{\alpha_r} g(y) |}{|x-y|^{d}} + \displaystyle \sup_{\substack{x,y \in I \\x \neq y }} \dfrac{| g(x) - g(y) |}{|x-y|^{d}}\\
    &= Lip_{d}(T^{\alpha_r} g) + Lip_{d}(g)\\
    &\le Lip_{d}(f) + \max_{i \in \mathbb{N}_{N}} \left(\frac{\|\alpha_{i,r}\|_{d}}{a_{i}^{d}}\right) ( Lip_{d}(g) + Lip_{d}(b_r) + \|g\|_\infty+\|b_r\|_\infty) + Lip_{d}(g)\\
    &\le Lip_{d}(f) + \left( 1 + \max_{i \in \mathbb{N}_{N}} \left(\frac{\|\alpha_i\|_{d}}{a_{i}^{d}}\right) \right) Lip_{d}(g) + 2 \max_{i \in \mathbb{N}_{N}} \left(\frac{\|\alpha_i\|_{d}}{a_{i}^{d}}\right) \|b_r\|_{d} \\
    &\le Lip_{d}(f) + \left( 1 + \max_{i \in \mathbb{N}_{N}} \left(\frac{\|\alpha_i\|_{d}}{a_{i}^{d}}\right) \right) Lip_{d}(g) + 2 \max_{i \in \mathbb{N}_{N}} \left(\frac{\|\alpha_i\|_{d}}{a_{i}^{d}}\right) \|b\|_{d} \label{eq:eq1}\tag{i},
    \end{align*}
    where $\|\alpha_i\|_{d}:= \sup\limits_{r \in \mathbb{N}} \|\alpha_{i,r}\|_{d}$.
    Also, 
    \begin{align*}
        |(T^{\alpha_r} g - g)(x)| &= |(f-g)(x)|+|\alpha_{i,r}(Q_{i}(x))|\cdot |(g - b_{r})(Q_{i}(x))|\\
        &\le \|f-g\|_\infty + \|\alpha\|_\infty \|g-b_r\|_\infty\\
        &\le \|f-g\|_\infty + \|\alpha\|_\infty (\|g\|_\infty + \|b_r\|_\infty)\\
        &\le \|f-g\|_{d} + \|\alpha\|_\infty (\|g\|_{d} + \|b_r\|_{d}).
    \end{align*}
    Hence,
    \begin{align*}
    \|T^{\alpha_r} g - g\|_{\infty} &\le \|f-g\|_{d} + \|\alpha\|_\infty (\|g\|_{d} + \|b\|_{d}).\label{eq:eq2}\tag{ii}
    \end{align*}
    Combining (i) and (ii), we get that the bound of $\|T^{\alpha_r}g - g\|_{d}$ is independent of $r$. Using Proposition \ref{Thm-nonstationary}, $\exists$ a unique $f_{b,Lip_d}^{\alpha} \in Lip_{d,f}(I)$ such that $f_{b,Lip_d}^{\alpha} = \displaystyle \lim_{ r \rightarrow \infty} T^{\alpha_1}\ o\ T^{\alpha_2}\ o \dots o\ T^{\alpha_r} g$ for any $g \in Lip_{d,f}(I)$. This completes the proof of the theorem.
\end{enumerate}
\end{proof}
\begin{definition}
The function $f_{b,Lip_d}^{\alpha}$ is called a Lipschitz non-stationary $\alpha$-fractal function with respect to $f,\alpha,b$ and the partition $\Delta.$ 
\end{definition}
\begin{remark}
As each $ T^{\alpha_r} $ is a contraction, there exists a unique stationary $\alpha$-fractal function $f_{r}^{\alpha}$ such that $ T^{\alpha_r}(f_{r}^{\alpha}) = f_{r}^{\alpha} $ and it satisfies the functional equation:
$$f_{r}^{\alpha}(x) = F_{i,r}\ (\ Q_{i}(x) ,\ f_{r}^{\alpha}(Q_{i}(x))))\ \ \forall\ \ x \in I_i,$$
where $Q_{i}(x) := {l_{i}}^{-1}(x)$. That is,
$$ f_{r}^{\alpha}(x) = f(x) + \alpha_{i,r}(Q_{i}(x)). f_{r}^{\alpha}(Q_{i}(x))- \alpha_{i,r}(Q_{i}(x)) b_{r}(Q_{i}(x)).$$
\end{remark}

\section{A Nonlinear Fractal Operator on $Lip_{d}(I)$}
Suppose $L_r: Lip_{d}(I) \to Lip_{d}(I)$ is a sequence of operators such that $\|L\|_\infty:= \displaystyle \sup_{r \in \mathbb{N}} \|L_r\|_\infty < \infty$ and satisfy $(L_r(f))(x_0) = f(x_0)$ and $(L_r(f))(x_N) = f(x_N)$. We set $b_r = L_{r}f$. The corresponding non-stationary $\alpha$-fractal function will be denoted by $f_{b}^{\alpha}$.

\begin{definition}
    Let $f \in Lip_{d}(I)$ and $\Delta$ be fixed. We define the $\alpha$-fractal operator $\mathfrak{F}_{b}^{\alpha} \equiv \mathfrak{F}_{\Delta,b}^{\alpha}$ as 
\begin{center}
$\mathfrak{F}_{b}^{\alpha} : Lip_{d}(I) \subset \mathcal{C}(I) \to \mathcal{C}(I) ,\ \ \mathfrak{F}_{b}^{\alpha}(f) = f_{b}^{\alpha}.$
\end{center}
\end{definition}
\begin{remark}
    In the case of a stationary fractal function, a similar construction is well studied in the literature \cite{Navascues1, Navascues2}. If we take $\alpha_{i,r} = \alpha_{i} \forall \;r \in \mathbb{N}, i \in \mathbb{N}_N$ and $b_r = Lf \;\forall \; r \in \mathbb{N}$, where $L: Lip_{d}(I) \to Lip_{d}(I)$ is an operator such that $(L(f))(x_0) = f(x_0)$ and $(L(f))(x_N) = f(x_N)$. Then the non-stationary $\alpha$-fractal function will coincide with the stationary one.
\end{remark}

Our next concern is to study the error approximation in the non-stationary perturbation process. The error bound in the different fractal approximations is well-studied in the stationary case \cite{Navascues2}. 

\begin{proposition}
    Let $f_{b}^{\alpha}$ be the non-stationary FIF corresponding to the seed function $f \in Lip_{d}(I)$. Then we have the following error bound
    \begin{equation}\label{err-eq1}
        \|f_{b}^{\alpha} - f\|_\infty \le \frac{\|\alpha\|_\infty}{1-\|\alpha\|_\infty} \sup_{r \in \mathbb{N}}\{ \|f-L_{r}(f)\|_\infty\}.
    \end{equation}
\end{proposition}

\begin{proof}
   The proof is similar to that given in Theorem 4.1. of \cite{NavascuesVerma}.
\end{proof}

\begin{corollary}
     Let $f \in Lip_{d}(I)$ be the germ function and $f_{b}^{\alpha}$ be the corresponding non-stationary FIF. Then for any $j \in \mathbb{N}$, we have the following inequality
     \begin{equation*}
         \|f_{b}^{\alpha} - L_j(f)\|_\infty \le \frac{1}{1-\|\alpha\|_\infty} \sup_{r \in \mathbb{N}}\{ \|f-L_{r}(f)\|_\infty\}.
     \end{equation*}
\end{corollary}

\begin{proof}
Let $j \in \mathbb{N}$. Using inequality (\ref{err-eq1}), we get
    \begin{align*}
       \|f_{b}^{\alpha} - L_j(f)\|_\infty &=  \|f_{b}^{\alpha} - f + f - L_j(f)\|_\infty\\
       &\le \|f_{b}^{\alpha} - f\|_\infty +\| f - L_j(f)\|_\infty\\
       &\le \frac{\|\alpha\|_\infty}{1-\|\alpha\|_\infty} \sup_{r \in \mathbb{N}}\{ \|f-L_{r}(f)\|_\infty\} + \| f - L_j(f)\|_\infty\\
       &\le \frac{1}{1-\|\alpha\|_\infty} \sup_{r \in \mathbb{N}}\{ \|f-L_{r}(f)\|_\infty\}.
    \end{align*}
\end{proof}

Based on the same arguments used in \cite{NavascuesVerma}, we know if $L_r$ is linear, then $\mathfrak{F}_{b}^{\alpha}$ is a linear operator. In order to keep track of this, let us write it down in the next proposition:

\begin{proposition}
    The fractal operator $\mathfrak{F}_{b}^{\alpha}$ is a linear operator, provided that the sequence of operators $L_r: Lip_{d}(I) \to Lip_{d}(I)$ are linear for each $r \in \mathbb{N}$.
\end{proposition}

Unless otherwise specified, note that we do not assume that $L_r$ is linear. As a result, the fractal operator is typically nonlinear (not necessarily linear). With regard to the conventional setting of fractal operators spread throughout the literature \cite{Navascues1, Navascues2}, the present findings abandon the general assumption of linearity and boundedness of the map $L_r$. Consequently, the research presented here may uncover possible applications of the fractal operator within the theory of unbounded and nonlinear operators.

Let us now collect some standard definitions of operators of interest in nonlinear functional analysis and perturbation theory. Let $A, B$ be two normed linear spaces.

\begin{definition}
    If an operator $\mathcal{T} :  A \to B$ maps bounded sets to bounded sets, then it is said to be topologically bounded.
\end{definition}

\begin{definition}
    Let $\mathcal{T}_1 : D(\mathcal{T}_1) \subset A \to B$ and $\mathcal{T}_2 : D(\mathcal{T}_2) \subset A \to B$ be two operators such that $D(\mathcal{T}_2) \subset D(\mathcal{T}_1)$. If $\mathcal{T}_1, \mathcal{T}_2$ satisfy the following inequality
$$\|\mathcal{T}_1(u)\| \le t_1\|u\| + t_2\|\mathcal{T}_2(u)\|\;\; \forall\; u \in D(\mathcal{T}_2),$$
where $t_1$ and $t_2$ are some non-negative constants, then $\mathcal{T}_1$ is said to be relatively (norm) bounded with respect to $\mathcal{T}_2$ or simply $\mathcal{T}_2$-bounded. The $\mathcal{T}_2$-bound of $\mathcal{T}_1$ is defined as the infimum of all possible values of $t_2$ satisfying the aforementioned inequality.
\end{definition}

\begin{definition}
    An operator $\mathcal{T} : A \to B $ is said to be Lipschitz if
there exists a constant $q > 0$ such that
$$\|\mathcal{T}(u) - \mathcal{T}(v)\| \le q\|u - v\|\ \forall\; u,v \in A.$$
For a Lipschitz operator $\mathcal{T}:  A \to B$, the Lipschitz constant of $\mathcal{T}$ is denoted by $|\mathcal{T}|$.
\end{definition}

\begin{definition}
    Let $\mathcal{T}_1 : D(\mathcal{T}_1) \subset A \to B$ and $\mathcal{T}_2 : D(\mathcal{T}_2) \subset A \to B$ be two operators such that $D(\mathcal{T}_2) \subset D(\mathcal{T}_1)$. If $\mathcal{T}_1, \mathcal{T}_2$ satisfies the following inequality 
$$\|\mathcal{T}_1(u) - \mathcal{T}_1(v)\| \le M_1\|u - v\| + M_2\|\mathcal{T}_2(u) - \mathcal{T}_2(v)\|\ \forall u,v \in D(\mathcal{T}_1),$$
where $M_1$ and $M_2$ are non-negative constants, then we say that $\mathcal{T}_1$ is relatively Lipschitz with
respect to $\mathcal{T}_2$ or simply $\mathcal{T}_2$-Lipschitz. The infimum of all such values
of $M_2$ is called the $\mathcal{T}_2$-Lipschitz constant of $\mathcal{T}_1$.
\end{definition}

\begin{proposition}
    The non-stationary fractal operator $\mathfrak{F}_{b}^{\alpha}: Lip_{d}(I) \to\mathcal{C}(I)$ is continuous whenever $L_r: Lip_{d}(I) \to Lip_{d}(I)$ is continuous for each $r \in \mathbb{N}$.
\end{proposition}

\begin{proof}
Let $(f_n)_{n \in \mathbb{N}}$ be a convergent sequence in $Lip_{d}(I)$, converges to $f \in Lip_{d}(I)$.
We have,
    \begin{equation*}
f_{b}^{\alpha} (x) = f(x) + \lim_{r \rightarrow \infty} \sum_{j=1}^{r} \alpha_{i,1}(Q_{i}(x))\dots\alpha_{i,j}(Q_{i}^{j}(x)) (f-L_j f)(Q_{i}^{j}(x)).
\end{equation*}
Now,
    \begin{align*}
    &|(f_n)_{b}^{\alpha} (x) - f_{b}^{\alpha} (x)|\\
    &\le |f_n(x) - f(x)| + |\lim_{r \rightarrow \infty} \sum_{j=1}^{r} \alpha_{i,1}(Q_{i}(x))\dots\alpha_{i,j}(Q_{i}^{j}(x)) (f_n - f - L_j f_n + L_j f)(Q_{i}^{j}(x))|\\
    &\le \|f_n - f\|_\infty + \lim_{r \rightarrow \infty} \sum_{j=1}^{r} \|\alpha\|_{\infty}^{j} ( \|f_n - f\|_\infty + \|L_j f_n - L_j f\|_\infty).
    \end{align*}
    Since the inequality holds for all $x \in I$, we have
    $$\|(f_n)_{b}^{\alpha} - f_{b}^{\alpha}\| \le \|f_n - f\|_\infty + \lim_{r \rightarrow \infty} \sum_{j=1}^{r} \|\alpha\|_{\infty}^{j} ( \|f_n - f\|_\infty + \|L_j f_n - L_j f\|_\infty).$$
    As the sequence $(f_n)$ converges to $f$, we get our desired result using continuity of $L_j, j \in \mathbb{N}$.
\end{proof}

\begin{proposition}
    If for each $r \in \mathbb{N}$, the operators $L_r: Lip_{d}(I) \to Lip_{d}(I)$ is a Lipschitz operator with Lipschitz constant $|L_r|$, then the non-stationary fractal operator $\mathfrak{F}_{b}^{\alpha}: Lip_{d}(I) \to \mathcal{C}(I)$ is also a Lipschitz operator, and $|\mathfrak{F}_{b}^{\alpha}| \le \dfrac{1 + |L| \|\alpha\|_\infty}{1-\|\alpha\|_\infty}$, where $|L| :=\displaystyle \sup_{r \in \mathbb{N}} |L_r| < \infty$.
\end{proposition}

\begin{proof}
Let $f,g \in Lip_{d}(I)$. Then
\begin{align*}
f_{b}^{\alpha} (x) &= f(x) + \lim_{r \rightarrow \infty} \sum_{j=1}^{r} \alpha_{i,1}(Q_{i}(x))\dots\alpha_{i,j}(Q_{i}^{j}(x)) (f-L_j f)(Q_{i}^{j}(x)),\\
g_{b}^{\alpha} (x) &= g(x) + \lim_{r \rightarrow \infty} \sum_{j=1}^{r} \alpha_{i,1}(Q_{i}(x))\dots\alpha_{i,j}(Q_{i}^{j}(x)) (g-L_j g)(Q_{i}^{j}(x)).
\end{align*}
Therefore,
\begin{align*}
    |f_{b}^{\alpha} (x) - g_{b}^{\alpha} (x)|&= |f(x) + \lim_{r \rightarrow \infty} \sum_{j=1}^{r} \alpha_{i,1}(Q_{i}(x))\dots\alpha_{i,j}(Q_{i}^{j}(x)) (f-L_j f)(Q_{i}^{j}(x))\\
    &\hspace{1.25cm}-g(x) - \lim_{r \rightarrow \infty} \sum_{j=1}^{r} \alpha_{i,1}(Q_{i}(x))\dots\alpha_{i,j}(Q_{i}^{j}(x)) (g-L_j g)(Q_{i}^{j}(x))|\\
    &\le |f(x) - g(x)| + |\lim_{r \rightarrow \infty} \sum_{j=1}^{r} \alpha_{i,1}(Q_{i}(x))\dots\alpha_{i,j}(Q_{i}^{j}(x)) (f - g - L_j f + L_j g)(Q_{i}^{j}(x))|\\
    &\le \|f - g\|_\infty + \lim_{r \rightarrow \infty} \sum_{j=1}^{r} \|\alpha\|_{\infty}^{j} ( \|f - g\|_\infty + \|L_j f - L_j g\|_\infty)\\
     &\le \|f - g\|_\infty + \lim_{r \rightarrow \infty} \sum_{j=1}^{r} \|\alpha\|_{\infty}^{j} ( \|f - g\|_\infty + |L_j|\cdot\| f - g\|_\infty)\\
     &\le \left(1 +  \sum_{j=1}^{\infty} \|\alpha\|_{\infty}^{j} ( 1 + |L|)\right)\cdot\| f - g\|_\infty\\
     &= \left(1 +  \dfrac{\|\alpha\|_{\infty}}{1 - \|\alpha\|_{\infty}} ( 1 + |L|)\right)\cdot\| f - g\|_\infty\\
     &= \dfrac{1 + |L|\cdot \|\alpha\|_\infty}{1-\|\alpha\|_\infty} \cdot\| f - g\|_\infty.
    \end{align*}
This holds for every $x \in I$, hence 
\[\|\mathfrak{F}_{b}^{\alpha}(f) - \mathfrak{F}_{b}^{\alpha}(g)\| = \|f_{b}^{\alpha} - g_{b}^{\alpha}\| \le \dfrac{1 + |L|\cdot \|\alpha\|_\infty}{1-\|\alpha\|_\infty} \| f - g\|_\infty.\]
This concludes the proof.
\end{proof}

\begin{proposition}
    The non-stationary fractal operator $\mathfrak{F}_{b}^{\alpha}: Lip_{d}(I) \to \mathcal{C}(I)$ is topologically bounded provided that $L_r: Lip_{d}(I) \to Lip_{d}(I)$ is uniformly bounded.
\end{proposition}    

\begin{proof}
Let $f$ be a function in $Lip_d(I)$.
    We have,
    \begin{equation*}
    \allowdisplaybreaks
    \begin{aligned}
        |f_{b}^{\alpha} (x)| &\leq  |f(x)| + |\lim_{r \rightarrow \infty} \sum_{j=1}^{r} \alpha_{i,1}(Q_{i}(x))\dots\alpha_{i,j}(Q_{i}^{j}(x)) (f-L_jf)(Q_{i}^{j}(x))|\\
        &\le \|f\|_\infty + \lim_{r \rightarrow \infty} \sum_{j=1}^{r} \|\alpha\|_{\infty}^{j} \|f-L_j f\|_\infty\\
        &\le \|f\|_\infty + \lim_{r \rightarrow \infty} \sum_{j=1}^{r} \|\alpha\|_{\infty}^{j} (\|f\|_\infty+\|L_j f\|_\infty)\\
        &\le (1 + \sum_{j=1}^{\infty} \|\alpha\|_{\infty}^{j}) \|f\|_\infty + \sum_{j=1}^{\infty} \|\alpha\|_{\infty}^{j} \|L_j f\|_\infty\\
        &= \dfrac{1}{1 - \|\alpha\|_{\infty}} \|f\|_\infty\ + \sum_{j=1}^{\infty} \|\alpha\|_{\infty}^{j} \|L_j f\|_\infty.
    \end{aligned}
    \end{equation*}
    Hence 
    \begin{equation*}
        \|\mathfrak{F}_{b}^{\alpha}(f)\|_\infty\ =\ \|f_{b}^{\alpha}\|_\infty\ \le\ \dfrac{1}{1 - \|\alpha\|_{\infty}} \|f\|_\infty\ + \sum_{j=1}^{\infty} \|\alpha\|_{\infty}^{j} \|L_j f\|_\infty.
    \end{equation*}
    Since $L_j (j \in \mathbb{N})$ is uniformly bounded, it follows from the above inequality that the operator $\mathfrak{F}_{b}^{\alpha}$ is topologically bounded.
\end{proof}

In the following propositions of this section, we assume that $L_r$ be a sequence of linear operators such that there exists a linear operator $L$ satisfying $\|L f\|_\infty = \sup\limits_{r \in \mathbb{N}} \|L_r f\|_\infty$. Let's move on to the following proposition using this presumption.

\begin{proposition}
     The non-stationary fractal operator $\mathfrak{F}_{b}^{\alpha}: Lip_{d}(I) \to \mathcal{C}(I)$ is relatively Lipschitz with respect to $L$ with $L$-Lispchitz constant of $\mathfrak{F}_{b}^{\alpha}$ not exceeding $\dfrac{\|\alpha\|_{\infty}}{1 - \|\alpha\|_{\infty}}$.
\end{proposition}

\begin{proof}
Let $f,g\in Lip_d(I)$. Then the functions will satisfy the following equation
    \begin{equation*}
f_{b}^{\alpha} (x) = f(x) + \lim_{r \rightarrow \infty} \sum_{j=1}^{r} \alpha_{i,1}(Q_{i}(x))\dots\alpha_{i,j}(Q_{i}^{j}(x)) (f-L_jf)(Q_{i}^{j}(x)).
\end{equation*}
Now,
\allowdisplaybreaks
\begin{align*}
    &|f_{b}^{\alpha} (x) - g_{b}^{\alpha} (x)|\\
    &= |f(x) + \lim_{r \rightarrow \infty} \sum_{j=1}^{r} \alpha_{i,1}(Q_{i}(x))\dots\alpha_{i,j}(Q_{i}^{j}(x)) (f-L_j f)(Q_{i}^{j}(x))\\
    &\hspace{1.25cm}- g(x) - \lim_{r \rightarrow \infty} \sum_{j=1}^{r} \alpha_{i,1}(Q_{i}(x))\dots\alpha_{i,j}(Q_{i}^{j}(x)) (g-L_j g)(Q_{i}^{j}(x))|\\
    &\le |f(x) - g(x)| + |\lim_{r \rightarrow \infty} \sum_{j=1}^{r} \alpha_{i,1}(Q_{i}(x))\dots\alpha_{i,j}(Q_{i}^{j}(x)) (f - g - L_j f + L_j g)(Q_{i}^{j}(x))|\\
    &\le \|f - g\|_\infty + \lim_{r \rightarrow \infty} \sum_{j=1}^{r} \|\alpha\|_{\infty}^{j} ( \|f - g\|_\infty + \|L_j f - L_j g\|_\infty)\\
    &\le \|f - g\|_\infty + \lim_{r \rightarrow \infty} \sum_{j=1}^{r} \|\alpha\|_{\infty}^{j} ( \|f - g\|_\infty + \|L_j(f - g)\|_\infty )\\
    &\le \left( 1 + \sum_{j=1}^{\infty} \|\alpha\|_{\infty}^{j} \right) \|f - g\|_\infty + \left( \sum_{j=1}^{\infty} \|\alpha\|_{\infty}^{j} \right) \|L(f - g)\|_\infty \\
    &= \left( \dfrac{1}{1 - \|\alpha\|_{\infty}} \right) \|f - g\|_\infty + \left( \dfrac{\|\alpha\|_{\infty}}{1 - \|\alpha\|_{\infty}} \right) \|Lf - Lg\|_\infty.
\end{align*}
For all $x$, the abovementioned inequality is true, hence
\begin{equation*}
    \|\mathfrak{F}_{b}^{\alpha}(f) - \mathfrak{F}_{b}^{\alpha}(g)\|_\infty = \|f_{b}^{\alpha} - g_{b}^{\alpha}\|_\infty \le \left( \dfrac{1}{1 - \|\alpha\|_{\infty}} \right) \|f - g\|_\infty + \left( \dfrac{\|\alpha\|_{\infty}}{1 - \|\alpha\|_{\infty}} \right) \|Lf - Lg\|_\infty.
\end{equation*}
This completes the proof.
\end{proof}

\begin{proposition}\label{relative bdd-proposition}
    The non-stationary fractal operator $\mathfrak{F}_{b}^{\alpha}: Lip_{d}(I) \to \mathcal{C}(I)$ is relatively bounded with respect to $L$ with $L$-bound is less than or equal to $\dfrac{\|\alpha\|_{\infty}}{1 - \|\alpha\|_{\infty}}$.
\end{proposition}

\begin{proof}
Let $f$ be an arbitrary function in $Lip_d (I)$.
    From (\ref{RB operator eq}), we have
    $$f_{b}^{\alpha} (x) = f(x) + \lim_{r \rightarrow \infty} \sum_{j=1}^{r} \alpha_{i,1}(Q_{i}(x))\dots\alpha_{i,j}(Q_{i}^{j}(x)) (f-L_jf)(Q_{i}^{j}(x)).$$
    Therefore,
    \begin{align*}
        |f_{b}^{\alpha} (x)| &= |f(x)| + |\lim_{r \rightarrow \infty} \sum_{j=1}^{r} \alpha_{i,1}(Q_{i}(x))\dots\alpha_{i,j}(Q_{i}^{j}(x)) (f-L_jf)(Q_{i}^{j}(x))|\\
        &\le \|f\|_\infty + \lim_{r \rightarrow \infty} \sum_{j=1}^{r} \|\alpha\|_{\infty}^{j} \|f-L_j f\|_\infty\\
        &\le \|f\|_\infty + \lim_{r \rightarrow \infty} \sum_{j=1}^{r} \|\alpha\|_{\infty}^{j} (\|f\|_\infty+\|L_j f\|_\infty)\\
        &\le \|f\|_\infty + \lim_{r \rightarrow \infty} \sum_{j=1}^{r} \|\alpha\|_{\infty}^{j} (\|f\|_\infty+\|L f\|_\infty)\\
        &= \|f\|_\infty + \dfrac{\|\alpha\|_{\infty}}{1 - \|\alpha\|_{\infty}} (\|f\|_\infty+\|L f\|_\infty)\\
        &= \dfrac{1}{1 - \|\alpha\|_{\infty}} \|f\|_\infty + \dfrac{\|\alpha\|_{\infty}}{1 - \|\alpha\|_{\infty}} \|L f\|_\infty.
    \end{align*}
    The aforementioned inequality holds for all $x$, hence
\begin{equation}\label{relative bdd-eq1}
   \|\mathfrak{F}_{b}^{\alpha}(f)\|_\infty\ =\ \|f_{b}^{\alpha}\|_\infty\ \le\ \dfrac{1}{1 - \|\alpha\|_{\infty}} \|f\|_\infty\ +\ \dfrac{\|\alpha\|_{\infty}}{1 - \|\alpha\|_{\infty}} \|L f\|_\infty.
\end{equation}
    This proves our claim.
    \end{proof}


\section{Stability and sensitivity analysis}
Let us now investigate the stability of the FIF with changeable parameters produced by IFS $ \mathcal{I}_r = \{ \mathcal{M}; W_{i,r}(x,y) = (l_{i}(x), F_{i,r}(x,y)), i \in \mathbb{N}_N \}$, where the maps are defined in (\ref{non-stationary IFS-eq1}) and $ \mathcal{M}=I\times [k_1,k_2]\subset \mathbb{R}^2$. The similar results for the stationary case can be observed in \cite{WangYu}. Let $\mathbf{\bar D}:= \{ (x_i, \bar{y}_i): i \in \mathbb{N}_{N}^{0} \}$ be another set of interpolation points in $\mathcal{M}$ which can be considered as the perturbations of ordinates of the points in $\mathbf{D}:= \{ (x_i, y_i)\in I\times [k_1,k_2]: i \in \mathbb{N}_{N}^{0} \}$. For the data set $\mathbf{\bar D}$, an IFS can be defined by
$ \mathcal{\bar I}_m = \{ \mathcal{M};\bar W_{i,r}(x,y) = (l_{i}(x), \bar F_{i,r}(x,y)), i \in \mathbb{N}_N \}$, where $l_{i}(x),\ i \in \mathbb{N}_N,$ are the maps defined in (\ref{non-stationary IFS-eq1}), and $\bar F_{i,r}$ are defined as
\begin{equation}\label{IFS-eq2}
    \bar F_{i,r} = \alpha_{i,r}(x) y + \hat{f}(l_{i}(x))-\alpha_{i,r}(x) \hat{b}_r(x),\ \ i \in \mathbb{N}_N,\ \ r \in \mathbb{N}.
\end{equation}
Here we consider the base functions $b_r$ and perturbed base functions $\bar b_r$ in $\mathcal{C}_f(I)$ such that $\sup\limits_{r\in \mathbb{N}}\Vert b_r\Vert_\infty<\infty$ and $\sup\limits_{r\in \mathbb{N}}\Vert \bar b_r\Vert_\infty<\infty$.
\begin{theorem}
    Let $\mathbf{D} := \{ (x_i , y_i) : i \in \mathbb{N}_{N}^{0} \}$ and $\mathbf{\bar D} := \{ (x_i , \bar{y}_i) : i \in \mathbb{N}_{N}^{0} \}$ be two data sets in $\mathcal{M}$. Let $f_{b}^{\alpha}$ be the non-stationary FIF for $\mathbf{D}$ generated by the sequence of IFSs $ \mathcal{I}_r = \{ \mathcal{M} ; W_{i,r}(x,y) = (l_{i}(x) , F_{i,r}(x,y)), i \in \mathbb{N}_N \}$ defined in (\ref{non-stationary IFS-eq1}) and $\bar f_{b}^{\alpha}$ be the non-stationary FIF for $\mathbf{\bar D}$ generated by the sequence of IFSs $ \mathcal{\bar I}_r = \{ \mathcal{M} ; \bar W_{i,r}(x,y) = (l_{i}(x) ,\bar F_{i,r}(x,y)), i \in \mathbb{N}_N \}$ defined through (\ref{IFS-eq2}). Then we have,
    \begin{equation}\label{sta-ineq}
        \|f_{b}^{\alpha} - \bar f_{b}^{\alpha}\|_\infty \le \dfrac{\|f - \hat{f}\|_\infty + \|\alpha\|_\infty \cdot \displaystyle \sup_{r \in \mathbb{N}} \{\|b_r - \hat{b}_r\|_\infty\}}{1 - \|\alpha\|_\infty}.
    \end{equation}
\end{theorem}
\begin{proof}
    From (\ref{self referential eq}), we have
    \begin{equation*}
f_{b}^{\alpha} (x) = f(x) + \lim_{r \rightarrow \infty} \sum_{j=1}^{r} \alpha_{i,1}(Q_{i}(x))\dots\alpha_{i,j}(Q_{i}^{j}(x)) (f-b_j)(Q_{i}^{j}(x)).
\end{equation*}
Therefore,
\begin{align*}
    &|f_{b}^{\alpha} (x) - \hat{f}_{b}^{\alpha} (x)|\\
    &= |f(x) + \lim_{r \rightarrow \infty} \sum_{j=1}^{r} \alpha_{i,1}(Q_{i}(x))\dots\alpha_{i,j}(Q_{i}^{j}(x)) (f-b_j)(Q_{i}^{j}(x)) - \hat{f}(x) \\
    & \hspace{1.25cm}- \lim_{r \rightarrow \infty} \sum_{j=1}^{r} \alpha_{i,1}(Q_{i}(x))\dots\alpha_{i,j}(Q_{i}^{j}(x)) (\hat{f}-\hat{b}_j)(Q_{i}^{j}(x))|\\
    &\le |f(x) - \hat{f}(x)| + |\lim_{r \rightarrow \infty} \sum_{j=1}^{r} \alpha_{i,1}(Q_{i}(x))\dots\alpha_{i,j}(Q_{i}^{j}(x)) (f - \hat{f} - b_j + \hat{b}_j)(Q_{i}^{j}(x))|\\
    &\le \|f - \hat{f}\|_\infty + \lim_{r \rightarrow \infty} \sum_{j=1}^{r} \|\alpha\|_{\infty}^{j} ( \|f - \hat{f}\|_\infty + \|b_j - \hat{b}_j\|_\infty)\\
    &\le \left( 1 + \sum_{j=1}^{\infty} \|\alpha\|_{\infty}^{j} \right) \|f - \hat{f}\|_\infty + \left( \lim_{r \rightarrow \infty} \sum_{j=1}^{r} \|\alpha\|_{\infty}^{j} \right) \sup_{r \in \mathbb{N}} \{\|b_r - \hat{b}_r\|_\infty\} \\
    &= \left( 1 + \dfrac{\|\alpha\|_{\infty}}{1 - \|\alpha\|_{\infty}} \right) \|f - \hat{f}\|_\infty + \left( \dfrac{\|\alpha\|_{\infty}}{1 - \|\alpha\|_{\infty}} \right) \sup_{r \in \mathbb{N}} \{\|b_r - \hat{b}_r\|_\infty\}\\
    &= \dfrac{\|f - \hat{f}\|_\infty + \|\alpha\|_\infty \cdot \displaystyle \sup_{r \in \mathbb{N}} \{\|b_r - \hat{b}_r\|_\infty\}}{1 - \|\alpha\|_\infty}.
\end{align*}
The above inequality holds for all $x \in I$; hence the inequality (\ref{sta-ineq}) follows.
\end{proof}
\begin{remark}
    Let $f,\hat{f}$ be two piecewise linear functions through the interpolation data sets $\mathbf{D}$ and $\mathbf{\bar{D}}$ respectively. Also assume that $ b_r=\hat{b}_r=b$ is a linear function passing through the points $(x_0,y_0)$ and $(x_N,y_N)$. Then we have
    $$\|f_{b}^{\alpha} - \bar f_{b}^{\alpha}\|_\infty \le \dfrac{1 + \|\alpha\|_\infty}{1 - \|\alpha\|_\infty} \max_{i \in \mathbb{N}_{N}^{0}} \{ |y_i - \bar y_i| \},$$
   which is the same result for stationary FIF given in \cite{WangYu}. So, our result can be treated as a generalisation of the existing result.
\end{remark}
\begin{remark}
    Perturbations of abscissas of interpolation points can be taken to affect the values of the non-stationary FIFs associated with the interpolation points. Also, perturbations of both abscissas and ordinates may be considered to examine the stability of the non-stationary FIF. For more details, the reader is invited to read the paper of Wang and Yu \cite{WangYu}.
\end{remark}
Next, we discuss the sensitivity of the non-stationary $\alpha$-FIF  defined by the IFS $\mathcal{I}_r$. Let $f,b_r,\alpha_{i,r}$ be as defined before and $T_{i,r}: \mathcal{M} \to \mathbb{R}, i \in \mathbb{N}_N, r \in \mathbb{N}$, be a sequence of continuous functions on $\mathcal{M}$ such that  for all $(x,y) \in \mathcal{M}$,
$$T_{i,r}(x,y) = f(x) + [\alpha_{i,r}(Q_{i}(x)) + t_{i,r} \theta_{i,r}(Q_{i}(x))] (g - b_{r})(Q_{i}(x)) + s_{i,r}\phi_{i,r}(Q_{i}(x)),$$
where $t_{i,r}, s_{i,r}$ are parameters of perturbation satisfying $|t_{i,r}|<1$ and $|s_{i,r}|<1$, $\phi_{i,r}, \theta_{i,r}\in Lip_d(I)$ satisfying $\max\limits_{i \in \mathbb{N}_N} \{\|\alpha_{i,r} + t_{i,r} \theta_{i,r}\|_\infty\} < 1$ and $\phi_{i,r}(x_0) = \phi_{i,r}(x_N) = 0.$ The function $T_{i,r}$ is a perturbation of the function $F_{i,r}$ for each $i \in \mathbb{N}_N, r \in \mathbb{N}.$ Thus the IFS $\mathcal{I'}_r = \{ \mathcal{M} ; (l_{i}(x) , T_{i,r}(x,y)), i \in \mathbb{N}_N \}$ may be treated as the perturbation IFS of the IFS $\mathcal{I}_r = \{ \mathcal{M} ; (l_{i}(x) , F_{i,r}(x,y)), i \in \mathbb{N}_N \}$. For each $ r \in \mathbb{N}, i \in \mathbb{N}_N, T_{i,r}$ is also contractive in the second variable and it satisfy
\begin{equation*}
    T_{i,r}(x_0,y_0) = y_{i-1}, \; \; T_{i,r}(x_N,y_N) = y_{i}.
\end{equation*}
Therefore the IFS $\mathcal{I'}_r = \{ \mathcal{M} ; (l_{i}(x) , T_{i,r}(x,y)), i \in \mathbb{N}_N \}$ determines a unique non-stationary FIF, denoted by $f_{b,s}^{\alpha,t}$.
\begin{proposition}
Let $\mathbf{D} := \{ (x_i , y_i) : i \in \mathbb{N}_{N}^{0} \}$ be a data set in $\mathcal{M}$. Let $f_{b}^{\alpha}$ be the non-stationary FIFs corresponding to the sequence of IFSs $ \mathcal{I}_r = \{ \mathcal{M} ; (l_{i}(x) , F_{i,r}(x,y)), i \in \mathbb{N}_N \}$ defined in (\ref{non-stationary IFS-eq1}) and  $f_{b,s}^{\alpha,t}$ be the non-stationary FIF generated by the sequence of IFSs $ \mathcal{ I'}_r = \{ \mathcal{M} ; (l_{i}(x) , T_{i,r}(x,y)), i \in \mathbb{N}_N \}$. Then
     \begin{equation} \label{sensitivity ineq}
         \|f_{b,s}^{\alpha,t}-f_{b}^{\alpha}\|_\infty \le \dfrac{\|\phi\|_\infty }{1-\|\alpha\|_\infty - \|t\|_\infty \|\theta\|_\infty}\|s\|_\infty + \dfrac{\|\theta\|_\infty \sup \{\|f-b_r\|_\infty\}}{(1-\|\alpha\|_\infty)(1-\|\alpha\|_\infty - \|t\|_\infty \|\theta\|_\infty)} \|t\|_\infty,
     \end{equation}
     where
     \begin{align*}
         &\|\phi\|_\infty = \sup_{r \in \mathbb{N}}\{ \max_{i \in \mathbb{N}_N} \|\phi_{i,r}\|_\infty\}, \; \; \|\theta\|_\infty = \sup_{r \in \mathbb{N}}\{ \max_{i \in \mathbb{N}_N} \|\theta_{i,r}\|_\infty\},\\
         &\|s\|_\infty = \sup_{r \in \mathbb{N}}\{ \max_{i \in \mathbb{N}_N} |s_{i,r}|\}, \; \; \|t\|_\infty = \sup_{r \in \mathbb{N}}\{ \max_{i \in \mathbb{N}_N} |t_{i,r}|\}.
     \end{align*}
\end{proposition}

\begin{proof}
    From (\ref{self referential eq}), we have
    \begin{equation} \label{self ref eq 1}
        f_{b}^{\alpha} (x) - f(x) = \lim_{r \rightarrow \infty} \sum_{j=1}^{r} \alpha_{i,1}(Q_{i}(x))\dots\alpha_{i,j}(Q_{i}^{j}(x)) (f-b_j)(Q_{i}^{j}(x)).
    \end{equation}
    For $r \in \mathbb{N}$, let us define an RB operator $V^{\alpha_r}$ on $\mathcal{M}$ by
    \begin{align*}
        (V^{\alpha_r} g)(x) =&~ T_{i,r}( {l_{i}}^{-1}(x) ,~g({l_{i}}^{-1}(x)))\\ =&~ f(x) + [\alpha_{i,r}(Q_{i}(x)) + t_{i,r} \theta_{i,r}(Q_{i}(x))] (g - b_{r})(Q_{i}(x)) + s_{i,r}\phi_{i,r}(Q_{i}(x))\\
        =&~ f(x) + [\mathbf{a}_r(x) + \mathbf{c}_r(x)] (g - b_{r})(Q_{i}(x)) + s_{i,r}\phi_{i,r}(Q_{i}(x)),
    \end{align*}
    where $\mathbf{a}_r(x) = \alpha_{i,r}(Q_{i}(x))$ and $ \mathbf{c}_r(x) = t_{i,r} \theta_{i,r}(Q_{i}(x))$.
    \begin{align*}
        &V^{\alpha_1} \circ V^{\alpha_2} \circ \dots \circ V^{\alpha_r}f(x)-f(x)\\&= \Big[\mathbf{a}_1(x) + \mathbf{c}_1(x)\Big]\Big(V^{\alpha_2} \circ V^{\alpha_3} \circ \dots \circ V^{\alpha_r}f- b_1\Big)(Q_{i}(x)) + s_{i,1} \phi_{i,1}(Q_{i}(x))
    \end{align*}
    Using induction, we obtain
    \begin{align*}
        &V^{\alpha_1} \circ V^{\alpha_2} \circ \dots \circ V^{\alpha_r}f(x)-f(x)\\
        &= \sum\limits_{j=1}^{r} \Big[\mathbf{a}_1(x) + \mathbf{c}_1(x)\Big]\Big[\mathbf{a}_2(x) + \mathbf{c}_2(x)\Big]\dots \Big[\mathbf{a}_j(x) + \mathbf{c}_{j}(x)\Big] (f-b_j) (Q_{i}^{j}(x))\\
        &\hspace{0.5cm}+\sum\limits_{j=1}^{r} s_{i,j} \phi_{i,j}(Q_{i}^{j}(x)) \Big[\mathbf{a}_1(x) + \mathbf{c}_1(x)\Big]\Big[\mathbf{a}_2(x) + \mathbf{c}_2(x)\Big] \dots\Big[\mathbf{a}_{j-1}(x) + \mathbf{c}_{j-1}(x)\Big],
    \end{align*}
    where $ Q_{i}^j$ is a suitable finite composition of mappings $Q_{i}.$ Now, taking the limit as $r \rightarrow \infty$, we get
    \begin{equation} \label{self ref eq 2}
    \begin{aligned}
        &f_{b,s}^{\alpha,t}(x) - f(x)\\
        &= \lim_{r \rightarrow \infty} \sum\limits_{j=1}^{r} \Big[\mathbf{a}_1(x) + \mathbf{c}_1(x)\Big]\Big[\mathbf{a}_2(x) + \mathbf{c}_2(x)\Big] \dots\Big[\mathbf{a}_j(x) + \mathbf{c}_{j}(x)\Big] (f-b_j) (Q_{i}^{j}(x))\\
        &\hspace{0.5cm}+ \lim_{r \rightarrow \infty} \sum\limits_{j=1}^{r} s_{i,j} \phi_{i,j}(Q_{i}^{j}(x)) \Big[\mathbf{a}_1(x) + \mathbf{c}_1(x)\Big] \Big[\mathbf{a}_2(x) + \mathbf{c}_2(x)\Big] \dots\Big[\mathbf{a}_{j-1}(x) + \mathbf{c}_{j-1}(x)\Big].
    \end{aligned}
    \end{equation}
    Subtracting (\ref{self ref eq 1}) from (\ref{self ref eq 2}), we get
    \begin{equation*}
    \begin{aligned}
        &f_{b,s}^{\alpha,t}(x) - f_{b}^{\alpha}(x)\\
         &= \lim_{r \rightarrow \infty} \sum\limits_{j=1}^{r} s_{i,j} \phi_{i,j}(Q_{i}^{j}(x)) \Big[\mathbf{a}_1(x) + \mathbf{c}_1(x)\Big] \Big[\mathbf{a}_2(x) + \mathbf{c}_2(x)\Big] \dots\Big[\mathbf{a}_{j-1}(x) + \mathbf{c}_{j-1}(x)\Big]\\
        &\hspace{0.5cm}+ \lim_{r \rightarrow \infty} \sum\limits_{j=1}^{r} \bigg[ \Big[\mathbf{a}_1(x) + \mathbf{c}_1(x)\Big]\Big[\mathbf{a}_2(x) + \mathbf{c}_2(x)\Big]\dots \Big[\mathbf{a}_j(x) + \mathbf{c}_{j}(x)\Big] \\
        &\hspace{2.4cm}-\mathbf{a}_1(x) \mathbf{a}_2(x)\dots\mathbf{a}_j(x) \bigg] (f-b_j) (Q_{i}^{j}(x))\\
        &= \lim_{r \rightarrow \infty} \sum\limits_{j=1}^{r} s_{i,j} \phi_{i,j}(Q_{i}^{j}(x)) \Big[\mathbf{a}_1(x) + \mathbf{c}_1(x)\Big] \Big[\mathbf{a}_2(x) + \mathbf{c}_2(x)\Big] \dots\Big[\mathbf{a}_{j-1}(x) + \mathbf{c}_{j-1}(x)\Big]\\
        &\hspace{0.5cm}+ \lim_{r \rightarrow \infty} \sum\limits_{j=1}^{r} \bigg[ \mathbf{a}_1(x)\cdot \mathbf{a}_2(x)\dots \mathbf{a}_{j-1}(x)\cdot \mathbf{c}_{j}(x) + \mathbf{a}_1(x)\cdot \mathbf{a}_2(x)\dots \mathbf{a}_{j-2}(x)\cdot \mathbf{c}_{j-1}(x)\\
        &\hspace{6cm}\times[\mathbf{a}_j(x)+ \mathbf{c}_{j}(x)]\\
        &\hspace{0.5cm}+ \mathbf{a}_1(x)\cdot \mathbf{a}_2(x)\dots \mathbf{a}_{j-3}(x)\cdot \mathbf{c}_{j-2}(x) [\mathbf{a}_{j-1}(x) + \mathbf{c}_{j-1}(x)] [\mathbf{a}_j(x) + \mathbf{c}_{j}(x)]\\
        &\hspace{0.5cm}+ \dots + \mathbf{c}_1(x) [\mathbf{a}_2(x) + \mathbf{c}_2(x)] [\mathbf{a}_3(x) + \mathbf{c}_3(x)]\dots [\mathbf{a}_j(x) + \mathbf{c}_{j}(x)] (f-b_j)(Q_{i}^{j}(x))
        \bigg].
    \end{aligned}
    \end{equation*}
    $\text{Let}\ \mathbf{a} = \sup_{r \in \mathbb{N}} \{\|\mathbf{a}_r\|_\infty\} = \|\alpha\|_\infty \ \text{and}\ \mathbf{c} = \sup_{r \in \mathbb{N}} \{\|\mathbf{c}_r\|_\infty\} = \|t\|_\infty \|\theta\|_\infty$.
Therefore,
\allowdisplaybreaks
\begin{align*}
    &|f_{b,s}^{\alpha,t}(x) - f_{b}^{\alpha}(x)|\\ 
    &\le \lim_{r \rightarrow \infty} \sum\limits_{j=1}^{r} \|s\|_\infty \|\phi\|_\infty (\mathbf{a} + \mathbf{c})^{j-1} + \sup \{\|f-b_r\|_\infty\} \\
    &\hspace{0.4cm}\times \lim_{r \rightarrow \infty} \sum\limits_{j=1}^{r} \bigg[ \mathbf{a}^{j-1} \cdot \mathbf{c} + \mathbf{a}^{j-2} \cdot \mathbf{c}
        (\mathbf{a} + \mathbf{c}) +\mathbf{a}^{j-3} \cdot \mathbf{c}
        (\mathbf{a} + \mathbf{c})^2 + \dots + \mathbf{c}
        (\mathbf{a} + \mathbf{c})^{j-1}
        \bigg],\\
        &= \|s\|_\infty \|\phi\|_\infty \lim_{r \rightarrow \infty} \sum\limits_{j=1}^{r} (\mathbf{a} + \mathbf{c})^{j-1} + \sup \{\|f-b_r\|_\infty\}\\
        &\hspace{1cm}\times\mathbf{c} \lim_{r \rightarrow \infty} \sum\limits_{j=1}^{r} \Big[ \mathbf{a}^{j-1} + \mathbf{a}^{j-2} \cdot
        (\mathbf{a} + \mathbf{c})
        + \mathbf{a}^{j-3} \cdot 
        (\mathbf{a} + \mathbf{c})^2 + \dots + (\mathbf{a} + \mathbf{c})^{j-1}
        \Big]\\
        &= \|s\|_\infty \|\phi\|_\infty \dfrac{1}{1-\mathbf{a} - \mathbf{c}} + \sup \{\|f-b_r\|_\infty\} \\
        &\hspace{0.4cm} \times \mathbf{c} \sum\limits_{j=1}^{\infty} (\mathbf{a} + \mathbf{c})^{j-1} \bigg[ 1 + \left(\dfrac{\mathbf{a}}{\mathbf{a} + \mathbf{c}}\right)+ \left(\dfrac{\mathbf{a}}{\mathbf{a} + \mathbf{c}}\right)^2 + \dots + \left(\dfrac{\mathbf{a}}{\mathbf{a} + \mathbf{c}}\right)^{j-1}\bigg]\\
        & =\dfrac{\|\phi\|_\infty \|s\|_\infty}{1-\mathbf{a} - \mathbf{c}} + \sup \{\|f-b_r\|_\infty\} \mathbf{c} \sum\limits_{j=1}^{\infty} (\mathbf{a} + \mathbf{c})^{j-1} \times \left(\dfrac{1-\left(\dfrac{\mathbf{a}}{\mathbf{a} + \mathbf{c}}\right)^{j}}{1-\left(\dfrac{\mathbf{a}}{\mathbf{a} + \mathbf{c}}\right)}\right)\\
        & =\dfrac{\|\phi\|_\infty \|s\|_\infty}{1-\mathbf{a} - \mathbf{c}} + \sup \{\|f-b_r\|_\infty\}\sum\limits_{j=1}^{\infty} \Big((\mathbf{a} + \mathbf{c})^{j} - (\mathbf{a})^{j}\Big)\\
        &= \dfrac{\|\phi\|_\infty \|s\|_\infty}{1-\mathbf{a} - \mathbf{c}} + \sup \{\|f-b_r\|_\infty\} \left(\dfrac{\mathbf{a} + \mathbf{c}}{1-\mathbf{a} - \mathbf{c}} - \dfrac{\mathbf{a} }{1-\mathbf{a}}\right)\\
        &= \dfrac{\|\phi\|_\infty }{1-\mathbf{a} - \mathbf{c}}\|s\|_\infty + \dfrac{\|\theta\|_\infty \sup \{\|f-b_r\|_\infty\}}{(1-\mathbf{a})(1-\mathbf{a} - \mathbf{c})} \|t\|_\infty.
\end{align*}
The above inequality holds for each $x \in \mathcal{M}$, hence 
\begin{equation*}
\begin{aligned}
    \|f_{b,s}^{\alpha,t}-f_{b}^{\alpha}\|_\infty \le &~ \dfrac{\|\phi\|_\infty }{1-\|\alpha\|_\infty - \|t\|_\infty \|\theta\|_\infty}\|s\|_\infty \\
    & + \dfrac{\|\theta\|_\infty \sup \{\|f-b_r\|_\infty\}}{(1-\|\alpha\|_\infty)(1-\|\alpha\|_\infty - \|t\|_\infty \|\theta\|_\infty)} \|t\|_\infty.
\end{aligned}
\end{equation*}
\end{proof}

\section{Continuous dependence on parameters $b,\alpha,\; \text{and}\; \triangle.$}
In this section, we will investigate the continuous dependence of the non-stationary $\alpha$-fractal function on different IFS parameters. The reader can refer to \cite{VermaViswanathan} for the same study in the stationary case. We will start with the continuous dependence of $ f^{\alpha}_{\Delta,b}$ on the sequence of base functions $b := \{b_r\}$.
\begin{theorem}
   Let $f \in \mathcal{C}(I),$ and the partition $ \triangle,$  sequence of scale functions $ \alpha_r \in \mathcal{C}(I),\ r \in \mathbb{N} $ with $ \| \alpha \|_{ \infty} < 1 $ are fixed. Let $ A =\{ b_{r} \in \mathcal{C}(I): b_{r}(x)=f(x), ~~~~\forall~ x = x_0, x_N \}.$ Then, the map $ \mathcal{A}: A \rightarrow \mathcal{C}(I)$ defined by $$\mathcal{A}(b)= f^{\alpha}_{\triangle, b}$$ is Lipschitz continuous.
   \end{theorem}
    \begin{proof}
   From Section \ref{sec-2.4}, we obtain that $f^{\alpha}_{\triangle,b}$ is unique for a fixed sequence of scale function $\alpha_r$, a partition $\triangle $, and a suitable sequence of the base function $b_r \in  \mathcal{C}(I)$. Further, $f^{\alpha}_{\triangle,b}$ satisfies the functional equation:
     for all $x \in I_{i},~i  \in\mathbb{N}_N$, we have
             \begin{equation*}
         f^{\alpha}_{\Delta,b}(x)= f(x) + \lim_{r \rightarrow \infty} \sum_{j=1}^{r} \alpha_{i,1}(Q_{i}(x))\dots\alpha_{i,j}(Q_{i}^{j}(x)) (f-b_j)(Q_{i}^{j}(x)).
        \end{equation*}
    Let $b_r, c_r \in A,~ \text{for}~r \in \mathbb{N}.$ Then
      $$ \mathcal{A}(b)(x)= f^{\alpha}_{\Delta,b}(x)= f(x) + \lim_{r \rightarrow \infty} \sum_{j=1}^{r} \alpha_{i,1}(Q_{i}(x))\dots\alpha_{i,j}(Q_{i}^{j}(x)) (f-b_j)(Q_{i}^{j}(x)), $$
     and $$ \mathcal{A}(c)(x)= f^{\alpha}_{\Delta,c}(x)= f(x) + \lim_{r \rightarrow \infty} \sum_{j=1}^{r} \alpha_{i,1}(Q_{i}(x))\dots\alpha_{i,j}(Q_{i}^{j}(x)) (f-c_j)(Q_{i}^{j}(x)).$$
   On subtraction, we get for $x \in I_{i},$
     \begin{equation*}
     \begin{aligned} 
     \mathcal{A}(b)(x)-\mathcal{A}(c)(x)=&~ \lim_{r \rightarrow \infty} \sum_{j=1}^{r} \alpha_{i,1}(Q_{i}(x))\dots\alpha_{i,j}(Q_{i}^{j}(x)) (c_j - b_j)(Q_{i}^{j}(x)).       
         \end{aligned}
         \end{equation*}
Therefore,
\allowdisplaybreaks
            \begin{align*} 
             \big|\mathcal{A}(b)(x)-\mathcal{A}(c)(x) \big|
              &\le \lim_{r \rightarrow \infty} \sum_{j=1}^{r} \big| \alpha_{i,1}(Q_{i}(x))\dots\alpha_{i,j}(Q_{i}^{j}(x)) (c_j - b_j)(Q_{i}^{j}(x)) \big|\\
              &\le \sum_{j=1}^{\infty} \big\| \alpha \big\|_{\infty}^{j} \|c_j - b_j\|_\infty\\
              &= \dfrac{\|\alpha\|_\infty}{1-\|\alpha\|_\infty} \| b - c \|_\infty.
              \end{align*}
For all $x \in I$, the aforementioned inequality holds. Hence, 
$$\|\mathcal{A}(b)-\mathcal{A}(c)\|_{\infty}
 \le \dfrac{\|\alpha\|_\infty}{1-\|\alpha\|_\infty} \| b - c \|_\infty.$$ 
 This shows that $\mathcal{A}$ is a Lipschitz continuous map with Lipschitz constant $ \frac{\|\alpha \|_{\infty} }{ 1-\|\alpha\|_{\infty}}.$ 
\end{proof}

\begin{theorem}
    Let $f \in \mathcal{C}(I),$ sequence of base function $b_r$ and the partition $\triangle,$ are fixed. Let $B=\{\alpha = \{\alpha_{r}\} : \alpha_{i,r} \in Lip_d(I)~\text{and}~ \| \alpha\|_{\infty} \le s   < 1 ~ \text{, where $s$ is a fixed number} \}$. Then the map $\mathcal{B}: B \to \mathcal{C}(I),$ defined by \[\mathcal{B}(\alpha)= f^{\alpha}_{\triangle, b}\] is continuous.
\end{theorem}
\begin{proof}
   For a fixed partition $\triangle $, a scale function $\alpha$, and a suitable sequence of base function $b_r$, the map $f^{\alpha}_{\triangle,b}$ is unique. Further, $f^{\alpha}_{\triangle,b}$ satisfies the functional equation: for all $x \in I_{i},~i  \in\mathbb{N}_N$ we have
   \begin{equation}
         f^{\alpha}_{\Delta,b}(x)= f(x) + \lim_{r \rightarrow \infty} \sum_{j=1}^{r} \alpha_{i,1}(Q_{i}(x))\dots\alpha_{i,j}(Q_{i}^{j}(x)) (f-b_j)(Q_{i}^{j}(x)).
    \end{equation}
    Let $ \alpha, \beta  \in B $, then from the above functional equation, we have
    $$ \mathcal{B}(\alpha)(x)= f(x) + \lim_{r \rightarrow \infty} \sum_{j=1}^{r} \alpha_{i,1}(Q_{i}(x))\dots\alpha_{i,j}(Q_{i}^{j}(x)) (f-b_j)(Q_{i}^{j}(x)),  $$
    and $$ \mathcal{B}(\beta)(x)= f(x) + \lim_{r \rightarrow \infty} \sum_{j=1}^{r} \beta_{i,1}(Q_{i}(x))...\beta_{i,j}(Q_{i}^{j}(x)) (f-b_j)(Q_{i}^{j}(x)).$$
    To show that $\mathcal{B}$ is continuous at $\alpha$, we subtract one from other of the above two equations, for $x \in I_{i},~i  \in\mathbb{N}_N,$ we have
    \begin{equation*}
    \begin{aligned} 
    &\mathcal{B}(\alpha)(x)-\mathcal{B}(\beta)(x)\\=& \lim_{r \rightarrow \infty} \sum_{j=1}^{r} \alpha_{i,1}(Q_{i}(x))\dots\alpha_{i,j}(Q_{i}^{j}(x)) (f-b_j)(Q_{i}^{j}(x))\\&~-\lim_{r \rightarrow \infty} \sum_{j=1}^{r} \beta_{i,1}(Q_{i}(x))\dots\beta_{i,j}(Q_{i}^{j}(x)) (f-b_j)(Q_{i}^{j}(x))\\
    =&\lim_{r \rightarrow \infty} \sum_{j=1}^{r} \Big(\alpha_{i,1}(Q_{i}(x))\dots\alpha_{i,j}(Q_{i}^{j}(x)) - \beta_{i,1}(Q_{i}(x))\dots\beta_{i,j}(Q_{i}^{j}(x)) \Big) (f-b_j)(Q_{i}^{j}(x))\\
    =& \lim_{r \rightarrow \infty} \sum_{j=1}^{r} \bigg( (\alpha_{i,1}-\beta_{i,1})(Q_i(x)) \beta_{i,2}(Q_{i}^{2}(x))\dots\beta_{i,j}(Q_{i}^{j}(x))\\
    &+ \alpha_{i,1}(Q_{i}(x)) (\alpha_{i,2}-\beta_{i,2})(Q_{i}^{2}(x)) \beta_{i,3}(Q_{i}^{3}(x))\dots\beta_{i,j}(Q_{i}^{j}(x))\\
    &+ \dots + \alpha_{i,1}(Q_{i}(x)) \alpha_{i,2}(Q_{i}^{2}(x))\dots\alpha_{i,j-1}(Q_{i}^{j-1}(x))(\alpha_{i,j}- \beta_{i,j})(Q_{i}^{j}(x))\bigg)(f-b_j)(Q_{i}^{j}(x)).
    \end{aligned}
    \end{equation*}
    Therefore,
\allowdisplaybreaks

  \begin{align*} 
  &\Big|\mathcal{B}(\alpha)(x)-\mathcal{B}(\beta)(x)\Big|\\
  \le & ~\lim_{r \rightarrow \infty} \sum_{j=1}^{r} \Big| (\alpha_{i,1}-\beta_{i,1})(Q_i(x)) \beta_{i,2}(Q_{i}^{2}(x))\dots\beta_{i,j}(Q_{i}^{j}(x))\Big|\\
  &+ \Big|\alpha_{i,1}(Q_{i}(x)) (\alpha_{i,2}-\beta_{i,2})(Q_{i}^{2}(x)) \beta_{i,3}(Q_{i}^{3}(x))\dots\beta_{i,j}(Q_{i}^{j}(x))\Big| + \dots\\
  &+ \Big|\alpha_{i,1}(Q_{i}(x)) \alpha_{i,2}(Q_{i}^{2}(x))\dots \alpha_{i,j-1}(Q_{i}^{j-1}(x))(\alpha_{i,j}- \beta_{i,j})(Q_{i}^{j}(x))\Big| \Big| (f-b_j)(Q_{i}^{j}(x))\Big|\\
  \le & \lim_{r \rightarrow \infty} \sum_{j=1}^{r} \Big(\|\alpha_{1}-\beta_{1}\|_\infty \|\beta_{i,2}\|_\infty \dots \|\beta_{i,j}\|_\infty
    + \|\alpha_{i,1}\|_\infty \|\alpha_{2}-\beta_{2}\|_\infty \|\beta_{i,3}\|_\infty \dots \|\beta_{i,j}\|_\infty+ \dots\\
    &+~ \|\alpha_{i,1}\|_\infty \|\alpha_{i,2}\|_\infty \dots \|\alpha_{i,j-1}\|_\infty \|\alpha_{j}-\beta_{j}\|_\infty\Big) \|f-b_j\|_\infty\\
    \le &~ \lim_{r \rightarrow \infty} \sum_{j=1}^{r}  \Big(\| \alpha - \beta\|_{\infty} \|\beta\|_{\infty}^{j-1} + \|\alpha\|_{\infty} \| \alpha - \beta\|_{\infty} \|\beta\|_{\infty}^{j-2} + \dots +\|\alpha\|_{\infty}^{j-1} \|\alpha - \beta\|_{\infty} \Big)\|f-b\|_\infty\\
    =& \|f-b\|_\infty \| \alpha - \beta\|_{\infty} \lim_{r \rightarrow \infty} \sum_{j=1}^{r} \|\alpha\|_{\infty}^{j-1} \bigg\{1 +\left(\dfrac{\|\beta\|_{\infty}}{\|\alpha\|_{\infty}}\right)+ \left(\dfrac{\|\beta\|_{\infty}}{\|\alpha\|_{\infty}}\right)^2 + \dots + \left(\dfrac{\|\beta\|_{\infty}}{\|\alpha\|_{\infty}}\right)^{j-1}\bigg\}.
    \end{align*}
          Without loss of generality, let $\|\beta\|_{\infty} < \|\alpha\|_{\infty}.$ Then
\begin{align*}
    &\Big|\mathcal{B}(\alpha)(x)-\mathcal{B}(\beta)(x)\Big|\\
    &\le \|f-b\|_\infty \| \alpha - \beta\|_{\infty} \sum_{j=1}^{\infty} \|\alpha\|_{\infty}^{j-1} \left(\dfrac{1-\left(\dfrac{\|\beta\|_\infty}{\|\alpha\|_\infty}\right)^{j}}{1-\left(\dfrac{\|\beta\|_\infty}{\|\alpha\|_\infty}\right)}\right)\\
    &=\|f-b\|_\infty \dfrac{\| \alpha - \beta\|_{\infty}}{\|\alpha\|_\infty - \|\beta\|_\infty} \sum_{j=1}^{\infty} \big( \|\alpha\|_{\infty}^{j} - \|\beta\|_{\infty}^{j}\big)\\
    &=\|f-b\|_\infty \dfrac{\| \alpha - \beta\|_{\infty}}{(1-\|\alpha\|_\infty)\cdot(1 - \|\beta\|_\infty)}\\
    &\le \| \alpha - \beta\|_{\infty} \dfrac{\|f-b\|_\infty}{(1-s)^2}.
\end{align*}
 The aforementioned inequality holds for all $x \in I$, therefore
 \begin{equation*}
     \Big|\mathcal{B}(\alpha)-\mathcal{B}(\beta)\Big| \le \| \alpha - \beta\|_{\infty} \dfrac{\|f-b\|_\infty}{(1-s)^2}.
 \end{equation*}
   Since $\|f^{\alpha}_{\triangle,b}\|_{\infty}$ is bounded and $\alpha $ is fixed, we have $\mathcal{B}$ is continuous at $\alpha.$ As $\alpha$ was taken arbitrarily, $\mathcal{B}$ is continuous on $B.$ 
    \end{proof}

Our next goal is to study the continuous dependence of $f^{\alpha}_{\Delta,b}$ on the partition $\Delta.$ 
In this regard, let us recall the following theorem.
\begin{theorem} \cite[Theorem 11.1]{Barnsley1} \label{refthm1}
Let $(X,d)$ be a complete metric space and $(P,d_p)$ be a metric space of parameters. Let $\{X;w_1,w_2,\dots w_N \}$ be a hyperbolic IFS with contractivity $c$. For $n\in \mathbb{N}_N$, let $w_n$ depend on the parameter $p \in (P,d_p)$ subject to the condition $d\big(w_{n_p}(x),w_{n_q}(x)\big) \le K ~d(p,q)$ for all $x \in X$ with $K$ independent of $n,p,$ or $x,$. Then the attractor $A(p) \in \mathcal{H}(X)$ depends continuously on the parameter $ p \in P$ with respect to the Hausdorff metric $h$ induced by $d$.
\end{theorem}

    

    \begin{remark} \label{partition}
Let $b_r:\mathcal{C}(I)\to \mathcal{C}(I)$ be such that $ \Vert b\Vert_\infty:=\sup\limits_{r\in \mathbb{N}}\Vert b_r\Vert_\infty<\infty$. Then from the bound of non-stationary $\alpha$-fractal function, we have $$\| f_{b}^{\alpha}\|_{\infty} \leq \|f\|_{\infty}+ \frac{\| \alpha\|_{\infty}}{1-\| \alpha\|_{\infty}}\|f-b\|_{\infty},$$
  where $\|f-b\|_{\infty}:= \displaystyle \sup_{r \in \mathbb{N}} \{\|f-b_r\|_{\infty}\}.$
  This demonstrates that the non-stationary $\alpha$-fractal function is constrained by a fixed value independent of partition $\triangle.$ For any acceptable partition of $I$, it is enough to work with $X = I \times [-R,R],$ where $R= \|f\|_{\infty}+ \frac{\| \alpha\|_{\infty}}{1-\| \alpha\|_{\infty}}\|f-b\|_{\infty}$.
  \end{remark}
    
    \begin{theorem}
    Let $f,~ b_r \in  Lip_{d}(I) $ be such that $ b_{r}(x_0)=f(x_0)$ and $ b_{r}(x_N)=f(x_N)$, with Lipschitz constants $k_f,~k_{b_r}$ respectively. Suppose that the scale function $\alpha_{i,r} \in Lip_{d}(I)$ with Lipschitz constant $k_{\alpha_r}$ such that $\|\alpha\|_{\infty}<1.$ Let the collection of all partitions of $I$ be denoted by $\mathcal{P}(I),$ that is,
    $$\mathcal{P}(I) = \{\triangle : \triangle = \{x_0,x_1,\dots,x_N\}~\text{such that}~ x_0<x_1<\dots<x_N\}.$$
    Then, the mapping $\mathcal{D}:\mathcal{P}(I) \to \mathcal{C}(I)$ defined by $\mathcal{D}(\triangle)=f^{\alpha}_{\triangle,b}$ is continuous.
    \end{theorem}
    \begin{proof}
    Let $\Delta:= \{x_i: i=0,1,2,\dots,N; x_0 < x_1< \dots < x_N\}$ and $\tilde{\Delta}:=\{\tilde{x}_i:i=0,1,2,\dots,N; x_0=\tilde{x}_0 < \tilde{x}_1< \dots < \tilde{x}_N=x_N\}$ be two partitions of $I$.
    Let $W_{i,r}(x,y)=\big(l_{i}(x),F_{i,r}(x,y)\big),$ where $l_i(x)$  and $F_{i,r}(x,y)$ as in (\ref{non-stationary IFS-eq1}) and $X=I \times [-R, R]$.  As the maps $l_i, F_{i,r}$ and $W_{i,r}$ depend on the partition chosen, hence we denote the above maps corresponding to the partition $\Delta$ by $l_{i}^{\Delta}, F_{i,r}^{\Delta}$ and $W_{i,r}^{\Delta}$ respectively.
    Therefore,  
    \begin{equation*}
    \begin{split}
    \big|l_{i}^{\Delta}(x)- l_{i}^{\tilde{\Delta}}(x)\big| =&~ \frac{1}{x_N-x_0} \big|(x_i-\tilde{x}_i)(x-x_0)+ (\tilde{x}_{i-1}-x_{i-1}) (x-x_N) \big|\\
    \le &~ |x_i-\tilde{x}_i| + |\tilde{x}_{i-1}-x_{i-1}|\\
    \le&~ 2 \| \Delta - \tilde{\Delta}\|_2,
    \end{split}
    \end{equation*}
    and
    \begin{equation*}
    \begin{split}
     \big|F_{i,r}^{\Delta}(x,y)-F_{i,r}^{\tilde{\Delta}}(x,y)\big|
     =&~\Big|\alpha_{i,r}(x)y+f\big(l_{i}^{\Delta}(x)\big)- \alpha_{i,r}(x)b_{r}(x)\\& - [\alpha_{i,r}(x)y+f\big(l_{i}^{\tilde{\Delta}}(x)\big)- \alpha_{i,r}(x)b_{r}(x)]\Big|\\
     \le&~  k_f \Big\|l_{i}^{\Delta}(x) - l_{i}^{\tilde{\Delta}}(x) \Big \|_2\\
    \le &~ 2 k_f \|\Delta - \tilde{\Delta}\|_2,
    \end{split}
    \end{equation*}
    We define a metric $\mu$ on $\mathbb{R}^2$ by
    $$\mu \big((x,y),(x',y')\big) = |x-x'|+\theta |y-y'|,$$
    where $\theta$ is a suitable constant mentioned below. 
    
    Using similar calculations as in Theorem 2.1. of \cite{Barnsley1}, we obtain that $\{X; W_{i,r}\}$ is a hyperbolic IFS with respect to the metric $\mu$ for 
    $$\theta < \dfrac{1-A}{ Rk_{\alpha}+A k_f  +\|\alpha \|_\infty k_b+ \|b\|_{\infty} k_{\alpha}},$$
    where $k_{\alpha} = \max \{ k_{\alpha_r} : r \in \mathbb{N}\},~ k_{b} = \max \{ k_{b_r} : r \in \mathbb{N}\}, A = \max \{a_i: i \in \mathbb{N} \}.$\\
Consequently,
    \begin{align*}
        \mu \big(W_{i,r}^{\Delta}(x,y),W_{i,r}^{ \tilde{\Delta}}(x,y)\big) &\le 2 \| \Delta - \tilde{\Delta}\|_2 + \theta 2 k_f \|\Delta - \tilde{\Delta}\|_2\\
        &= 2(1+\theta k_f) \|\Delta - \tilde{\Delta}\|_2.
    \end{align*}
    Therefore, the IFS maps $W_{i,r}$ depend continuously on the partition $ \Delta \in \mathcal{P}(I)$. Consequently, Theorem \ref{refthm1} asserts that the attractor $G(\triangle) \in \mathcal{H}(X)$  depends continuously on $ \Delta \in \mathcal{P}(I)$ with respect to the Hausdorff metric $h$. 
      Therefore, for given $ \epsilon > 0$ and $\Delta \in \mathcal{P}(I),$ we have  $$ \|   \mathcal{D}(\tilde{\Delta})-\mathcal{D} (\Delta)\|_\infty= \|f_{\tilde{\Delta},b}^\alpha-f_{\Delta,b}^\alpha\|_{\infty}  < \epsilon ~~\text{whenever}~~ \| \tilde{\Delta} -\Delta\|_2 < \delta.$$ That is, $\mathcal{D}$ is continuous at $\Delta.$ Since $\Delta \in \mathcal{P}(I)$ is arbitrary, $\mathcal{D}$ is continuous on $\mathcal{P}(I)$.
\end{proof}
\section*{Conclusion}
 In this article, we studied several analytical properties of the non-stationary fractal operator corresponding to the non-stationary FIFs. Note that in the construction of our proposed interpolant, the crucial IFS parameters
are the sequence of base functions and scaling functions. In literature, it is studied that the fractal dimension of FIFs (stationary) depends on the IFS parameters. Therefore, we attempt to compute the fractal dimensions of the non-stationary interpolant in our future investigation by selecting suitable parameters.
\bibliographystyle{amsplain}
\bibliography{references}

\end{document}